\documentclass[final]{siamart1116}
\usepackage{verbatim}
\usepackage{amsmath,bm}
\usepackage{graphicx}

\usepackage{epstopdf}
\ifpdf
  \DeclareGraphicsExtensions{.eps,.pdf,.png,.jpg}
\else
  \DeclareGraphicsExtensions{.eps}
\fi

\usepackage{amssymb}
\usepackage{amsmath,color}

\usepackage{comment}

\numberwithin{theorem}{section}

\newtheorem{rmk}{Remark}

\newcommand{\ba}{\begin{array}}
\newcommand{\ea}{\end{array}}
\newcommand{\bea}{\begin{eqnarray}}
\newcommand{\eea}{\end{eqnarray}}
\newcommand{\be}{\begin{equation}}
\newcommand{\ee}{\end{equation}}
\newcommand{\bd}{\begin{displaymath}}
\newcommand{\ed}{\end{displaymath}}
\newcommand{\bi}{\begin{itemize}}
\newcommand{\ei}{\end{itemize}}
\newcommand{\bn}{\begin{enumerate}}
\newcommand{\en}{\end{enumerate}}
\newcommand{\tbox}[1]{{\mbox{\tiny \rm #1}}}
\newcommand{\mbf}[1]{{\mathbf #1}}
\newcommand{\ci}{\cite}
\newcommand{\pa}{\partial}
\newcommand{\f}{\frac}

\newcommand{\dx}{h}               
\newcommand{\dt}{{\Delta t}}       
\newcommand{\bigO}{{\mathcal O}}       
\newcommand{\RR}{\mathbb{R}}
\newcommand{\nth}{{n_\theta}}      
\newcommand{\nph}{{n_\phi}}      
\newcommand{\Jni}{J_{\tbox{near},i}}
\newcommand{\bmu}{{\bm \mu}}        
\newcommand{\sfrac}[2]{{\textstyle\frac{#1}{#2}}}

\begin{document}

\newcommand{\TheTitle}{High-order discretization of a stable time-domain integral equation for 3D acoustic scattering} 
\newcommand{\TheAuthors}{A. H. Barnett, L. Greengard, and T. Hagstrom}  

\headers{Time-domain integral equation for 3D acoustic scattering}{\TheAuthors}

\title{{\TheTitle}\thanks{
\funding{Supported in part by NSF Grant DMS-1418871 and the U.S. Department
  of Energy ASCR Applied Mathematics program. Any conclusions or
  recommendations expressed in this paper are those of the authors 
and do not necessarily reflect the views of the NSF or DOE.}}}
\author{Alex Barnett%
  \thanks{Flatiron Institute, New York, NY
    (\email{abarnett@flatironinstitute.org}).}%
  \and
  Leslie Greengard%
  \thanks{Flatiron Institute and Courant Institute of Mathematical Sciences,
    New York, NY (\email{lgreengard@flatironinstitute.org}).}
  \and
  Thomas Hagstrom%
  \thanks{Department of Mathematics,
  Southern Methodist University, Dallas, TX (\email{thagstrom@smu.edu})}
}


\maketitle

\begin{abstract}
We develop a high-order, explicit method for acoustic scattering in 
three space dimensions based on a combined-field time-domain integral 
equation. The spatial discretization, of Nystr\"om type, uses Gaussian 
quadrature on panels combined with a special treatment of the
weakly singular kernels arising in near-neighbor interactions. 
In time, a new class of convolution splines is used in a 
predictor-corrector algorithm. Experiments on a torus and a perturbed 
torus are used to explore the stability and accuracy of the proposed 
scheme.
This involved around one thousand solver runs,
at up to 8th order and up to around 20,000 spatial unknowns,
demonstrating 5--9 digits of accuracy.
In addition we show that parameters in the combined field 
formulation, chosen on the basis of analysis for the sphere and other 
convex scatterers, work well in these cases.  
\end{abstract}



\section{Introduction}
\label{s:intro}

Problems involving the scattering of waves by obstacles have countless
applications in science and technology. For time-harmonic data, numerical methods
based on integral equations are well-developed and widely used.
Advantages of using an integral equation formulation,
rather than the partial differential equation itself, include 
superior conditioning when second-kind formulations are used,
reduced dimensionality, the lack of a need for mesh generation in the volume,
and the availability of fast solvers \ci{wideband3d,Chewfastbook,rokhlin_1993}.
They also impose outgoing radiation conditions exactly, 
avoiding the need for artificial boundary conditions on a truncated 
computational domain when considering exterior problems.

By contrast, numerical methods for time-domain integral equations are not 
so widely used, and the supporting theory is not nearly as complete. 
This is due, in large part, to the dependence of the relevant layer
potentials on their space-time history, making them somewhat unwieldy
in the absence of fast algorithms (see Remark \ref{r:fast} below).

Surprisingly, much of the literature is focused on equations
involving the single layer potential 
(see, for example, the review by Ha-Duong \ci{RetardRev} and the 
recent monograph by Sayas \ci{TDIE}), which lead to first kind equations 
in either the time or frequency
domains. Some exceptions that make use of second kind formulations 
include \cite{Ulku2013}, in which a quasi-explicit marching scheme is 
developed for the 
time domain magnetic field integral equation, 
and \cite{Valdes2013}, in which high-order accurate
Calderon-preconditioned versions of the electric field integral
equation are developed.

Here, we focus on the scalar wave equation; our main contributions are (1) 
the development of a high-order Nystr\"om discretizations for the combined field
integral equation proposed in \cite{TDIEstab} and (2) a numerical study of
the stability of explicit marching schemes for the resulting system.
In the frequency domain, combined field
equations are used to avoid singularities or near-singularities associated 
with eigenvalues
of the interior Helmholtz problem \cite{kress91} \cite[p.~48--49]{coltonkress}. For time-domain calculations, it is argued in
\cite{TDIEstab} that these interior resonances will always dominate the long time behavior
of the densities, rendering the recovery of the solution in the volume increasingly
ill-conditioned. Similar conclusions for electromagnetic problems are reached 
by Shanker et al in \ci{MichetalCFIE}, and potential benefits of direct formulations
to avoid numerical dispersion appear in the work of Banjai \ci{BanjaiCQ}.

More precisely, we consider the computation of
the scattered wave $u(x,t)$ inducing by an incoming acoustic wave impinging
on an obstacle in a uniform medium. The function $u$ satisfies
the equation
\be
\f {\pa^2 u}{\pa t^2} = c^2 \nabla^2 u, \qquad u(x,0)=\f {\pa u}{\pa t}(x,0)=0, \label{waveq}
\ee
for $(x,t) \in \Omega \times [0,T]$ where $\Omega \subset \mathbb{R}^3$ is the 
domain exterior to a compact obstacle with boundary $\Gamma$. 
In this work, for definiteness, we set unit sound speed $c=1$ and consider the 
Dirichlet (i.e.\ sound-soft) problem,
\bd
u(x,t)=g(x,t), \ \ \ x \in \Gamma , \; t \in [0,T], \label{Dbc}
\ed
where $g = - u_\tbox{inc}$ on the surface, and $u_\tbox{inc}$ is a given incident wave.
Our methods could also be applied to the Neumann (sound-hard) problem. 

For a target point $x \in \Omega$ and target time $t$, the retarded single and double layer potentials applied to
a density $\mu(y,t)$, $y \in \Gamma$, are defined by \cite[Sec.~10.7]{guentherlee} \cite[Ch.~1]{TDIE}
\begin{eqnarray}
\mathcal{S} \mu(x,t) & = & \int_{\Gamma}\frac{\mu(y,t-|x-y|)}{4\pi |x-y|}dS_y, \label{Slp} \\
\mathcal{D} \mu(x,t) & = & \int_{\Gamma}\frac{n_y\cdot(x-y)}{4\pi |x-y|}\left[
\frac{\mu(y,t-|x-y|)}{|x-y|^2}+\frac{\f {\pa \mu}{\pa t} (y,t-|x-y|)}{|x-y|}\right]dS_y .
\label{Dlp}
\end{eqnarray}
Introducing surface weight functions $a(y)$, $b(y)$, for $y\in\Gamma$
(analogous to the combined field parameter $\eta$ in the time-harmonic case \cite{kress91}),
we represent the scattered wave $u$ in the form
\be
u(x,t)\;=\;\mathcal{D} \mu(x,t) + \mathcal{S}\left(a \f {\pa \mu}{\pa t} + b \mu \right) (x,t) .
\label{solrep}
\ee
Taking the limit of (\ref{solrep}) as $x \rightarrow \Gamma$, using
the jump relation for the double-layer operator \cite[Sec.~1.3--1.4]{TDIE}, we obtain a
linear integral equation for $\mu(x,t)$, $x \in \Gamma$,
\be
\f {\mu(x,t)}{2} + D \mu (x,t) + S \left(a \f {\pa \mu}{\pa t} + b \mu \right) (x,t)
\;= \; g(x,t), \qquad x\in\Gamma, \; t\in[0,T],
\label{cftdie}
\ee
where $D$ is the principal value part of the double layer potential (\ref{Dlp}) spatially restricted to
$\Gamma \times \Gamma$, and $S$ is the
single layer potential (\ref{Slp}) spatially restricted to 
$\Gamma \times \Gamma$.
When $\Gamma$ is smooth, both $D$ and $S$ have weakly singular kernels with 
singularities bounded by $1/|x-y|$, as would be the case for the Laplace equation.
Note that \eqref{cftdie} is of Volterra type in time, that is, at time $t$
the $D$ and $S$ operators involve only the density history $\mu(\cdot,t')$ for $t'<t$. The first term suggests that it is also of the second kind.
However, even for the case $a\equiv b\equiv 0$, it does not
fall into the standard Fredholm theory \cite[Sec.~10.7]{guentherlee},
since $D$ is not compact.
Rigorous proofs of existence and uniqueness, along with some regularity estimates, for both a single and double layer formulation, may be based on the combination of Laplace transformation in time and subsequent analysis of the parametrized spatial integral equations \ci{RetardRev,TDIE}. 
Direct proofs of convergence can also be carried out (using the double layer alone) 
by fixed point iteration \cite{guentherlee}.

The choice of the parameters $a$ and $b$ and their effect on the long time behavior of
$\mu$ is a focus of \ci{TDIEstab}. There it is shown that:
\begin{description}
\item[i.] For convex obstacles, $a(y)\equiv 1$ is a natural choice, since an asymptotic analysis at
high frequency indicates that it leads to an optimal cancellation of the leading part of the
delay term; see additional discussions in Section \ref{TempD};
\item[ii.]
  Any positive $b$ is sufficient to remove the damaging zero-frequency
  interior Neumann resonance.
  For $b$ sufficiently large, the long time exponential decay rate of the density $\mu$ will match
that of the dominant physical scattering resonances, though for cases other than the sphere an
optimal choice of this parameter is unclear. For the sphere one should take $b\equiv R^{-1}$ where $R$ is its radius. 
\end{description}

Below, we examine the effect of simple choices for $a$ and $b$ for
scattering by a nonconvex obstacle. In particular, numerical experiments in
Section \ref{NumExp} compare results for various $a$ and $b$ when $\Gamma$ is 
the boundary of a torus or a perturbed torus.

Since the boundary is time-invariant, the integral equation (\ref{cftdie}) is
separable in space and time, so it is natural to separate the spatial
and temporal discretizations. Here we propose a high-order
Nystr\"om-like approximation in
space combined with a straightforward time marching procedure.
For this, we introduce a set of $N$ surface nodes, $x_j \in \Gamma$,
$j=1, \ldots ,N$ and a regular sequence of time steps $t_k = k \Delta t \in [0,T]$,
and represent the density $\mu$ by interpolation from its discrete values
\bd
\mu_j^k \;\approx\; \mu(x_j,t_k) ~.
\ed
To evolve forward in time,
the $k$th time-step consists of applying an explicit, fixed linear rule which gives the vector $\{\mu_j^k\}_{j=1,\ldots,N}$
in terms of the history $\{\mu_j^{k-r}\}_{j=1,\ldots,N}^{r=1,\ldots,n}$ and the current
data vector $\{g(x_j,t_k)\}_{j=1,\ldots,N}$.
Because of the strong Huygens principle, the number of previous times $n$ needed is
essentially the diameter of the obstacle divided by $\Delta t$.
The contributions of these past values come from approximating the retarded 
spatial integrals in (\ref{cftdie}).
We account for this history dependence by constructing
matrices $S^\dx$, $D^\dx$ and $W^\dx$ such that
for each target node $x_i$, 
\begin{eqnarray}
\int_{\Gamma} \frac{f(y)}{4\pi |x_i-y|}dS_y & = & \sum_{\ell=1}^{N'} S^\dx_{i\ell} f(y_\ell)
+ \bigO(\dx^{\gamma}) , \nonumber \\ 
\int_{\Gamma} \frac{n_y\cdot(x_i-y)}{4\pi |x_i-y|^3} f(y) dS_y & = & \sum_{\ell=1}^{N'} D^\dx_{i\ell}
f(y_\ell) + \bigO(\dx^{\gamma}), \label{matrices} \\
\int_{\Gamma} \frac{n_y\cdot(x_i-y)}{4\pi |x_i-y|^2} f(y) dS_y & = & \sum_{\ell=1}^{N'} W^\dx_{i\ell}
f(y_\ell) + \bigO(\dx^{\gamma}) , \nonumber
\end{eqnarray}
where $\dx = \bigO(N^{-1/2})$ is a measure of the spatial grid spacing,
and the convergence order $\gamma$ is determined by the scheme.
Here $\{y_\ell\}_{\ell=1,\ldots,N'}$ is a set of $N'>N$ nodes,
described in Section \ref{SpaceD},
comprising
all nodes $x_j$ far from $x_i$, plus a set of
auxiliary nodes designed to handle the
weakly singular contribution near to $x_i$ (see Fig.~\ref{f:geom}(c)).
The point samples $f(y_{\ell})$ in the expressions above involve the density and its time derivative at retarded times, e.g.
\be
f(y_{\ell}) \;=\; \mu(y_{\ell},t_k-|x_i-y_{\ell}|), \quad \f {\pa \mu}{\pa t}
(y_{\ell},t_k-|x_i-y_{\ell}|) ,
\label{f}
\ee
so that the retarded evaluation times do not correspond to points at previously
computed time steps, $t_{k-r}$.
To handle this, and to approximate the time derivatives, we introduce temporal interpolants
at each spatial grid point. As also suggested by Davies and Duncan
\ci{DaviesDuncanSpline,DaviesDuncanVolterraSpline}, our interpolants take the form 
\be
\mu_j (t_k-\tau) = \sum_{r} \omega_r (\tau) \mu_j^{k-r} ,
\label{ConvApp}
\ee
where $t_k$ is the current time.
However, we use different basis functions, $\omega_r$.
The temporal interpolating functions
and some properties of the time-stepping schemes defined 
in \ci{DaviesDuncanSpline,DaviesDuncanVolterraSpline}
are discussed in Section \ref{TempD}.  

The experiments in Sections \ref{SpaceD} and \ref{NumExp} verify the high-order convergence of the 
fully discrete algorithm.
In addition, the stability properties of the time marching scheme, in its
explicit predictor-corrector form, are investigated, revealing a 
rather surprising ``inverse CFL" constraint
of the form
\be
\dt \ge c_\tbox{iCFL} \dx ~.
\label{iCFL}
\ee
Such conditions have been noted before, but typically in the context of finite
volume methods \cite{bochev,burman,kuther}.

We note that most authors use fully implicit
time marching schemes, although there are exceptions \cite{Ulku2013}.
While the linear systems to be solved involve only local
interactions, a direct solver on a surface, even with optimal ordering, is likely to require
$\bigO(N^{3/2})$ work to factor,
with subsequent solves required at each time step requiring $\bigO(N \ln{N})$ flops.
Clearly, an explicit method (including predictor-corrector iterations) is cheaper, 
as no matrix factorizations are necessary
and the cost per time step is $\bigO(N)$. 
Finally, in Section \ref{Conclude} we summarize our results, and point to future 
enhancements and generalizations of our method. 

\begin{rmk}[Software]   
An open source MATLAB (with some Fortran90) implementation of the methods
of this paper is freely available from the repository 
\url{https://github.com/ahbarnett/BIE3D},
with code for most of the
figures from this paper in the {\tt timedomainwaveeqn/paper} directory.
\label{r:code}
\end{rmk}

\begin{rmk}[Fast algorithms] 
A critical bottleneck in using retarded layer potentials for the solution  
of the wave equation is that they require 
$\bigO(N^2)$ memory and $\bigO(N^2)$ work per timestep. Fortunately, over the
last two decades fast algorithms have been developed which reduce both of 
these costs to $\bigO(N \log N)$. These are described, for example, in 
\cite{Bleszynski,PWFTD,Chewfastbook,MengBoagMichielssen,Shanker2003,Yilmaz2004} 
and can be used to accelerate 
most Galerkin or Nystr\"om discretization schemes, including the method developed here.
In the present work, we make use of direct summation methods for the sake of simplicity.
\label{r:fast}
\end{rmk}

\section{Temporal discretization}\label{TempD}

Two commonly used approaches to time-stepping time-domain integral equations are Galerkin methods,
discussed extensively in Ha-Duong's review \ci{RetardRev}, and convolution quadrature \cite{lubich},
discussed extensively in Sayas' monograph \ci{TDIE}. Although provably stable with
exact integration, Galerkin methods have been found to be sensitive to quadrature errors arising
from the need to compute integrals in cut elements when discontinuous polynomial bases in
time are employed. As a result, various smooth nonstandard basis functions have been proposed
\ci{MOT1,SauterVeit}. Convolution quadrature methods, on the other hand, have been found to be
more robust. Their construction, based on combining standard time-marching schemes for
ordinary differential equations with representations of convolution in the Laplace domain,
leads to methods which do not respect the strong Huygens principle obeyed by the layer potentials
(\ref{Slp})--(\ref{Dlp}). That is, the solution updates involve the entire time history of the
potentials. As such, special methods must be introduced to alleviate the storage and
computation costs for long time computations; see, e.g., \ci{BanjaiCQ}. 

Most relevant to our approach is the convolution spline method of Davies and Duncan 
\ci{DaviesDuncanSpline,DaviesDuncanVolterraSpline}. They present a simple reinterpretation of
convolution quadrature as a temporal approximation (\ref{ConvApp}), and suggest taking
$\omega_r(t)$ to be smooth,
compactly-supported splines (both standard B-splines and more exotic bases), thus restoring 
the finite time history of the layer potential operators. However, their approximations
(\ref{ConvApp}) are quasi-interpolatory; that is, the basis functions $\omega_r$ do {\em not
satisfy}
\be
\omega_r (t_k) = \delta_{rk} . \label{Lagrange}
\ee
As a result, their proposed methods are limited to second order accuracy, even if
B-splines of high degree are used.
Higher order accuracy can be achieved using ``marching on in time" (MOT) schemes, which are more
closely related to the approach presented here (see \cite{Valdes2013} and references therein).
We propose the use of smooth 
piecewise polynomial bases which are designed to satisfy (\ref{Lagrange}). The functions, which
we term ``difference splines'' (or D-splines for short) \ci{Dspline} are defined as follows.

\begin{description}
\item[i.] Let $q$ be an integer and $\{\tau_k\}_{k=0,1,\dots}$ be a set of nodes, with $\tau_0=0$, and $\tau_{k+1}>\tau_k$ for each $k$. These nodes need not be the $t_k$ from the
previous section.
Let $S_k$ be the
set of $2q+1$ nearest nodes to be used in the difference stencil
for $\tau_k$. There are two cases:
\begin{description}
\item[a.] (Interior node, $k\ge q$): $S_k:=\{\tau_{k-q}, \ldots ,\tau_{k+q} \}$.
\item[b.] (Boundary node, $k<q$): $S_k:=\{\tau_0, \ldots , \tau_{2q}\}$.
\end{description}
\item[ii.] Let $P_k(\tau)$ be the degree $2q$ Lagrange interpolating polynomial defined by data, $\{ y_r \}$, 
on the stencil $S_k$.
That is, $P_k(\tau)$ is the unique polynomial of degree $2q$ satisfying $P_k(\tau_r)=y_r$ for all $r$ with $\tau_r \in S_k$.
\item[iii.] On the interval $(\tau_k,\tau_{k+1})$ the D-spline interpolant of the data $\{ y_r \}$,
$D_{k+1/2}(\tau)$, is defined as the degree $4q+1$
Hermite interpolant of $P_k(\tau)$ and $P_{k+1}(\tau)$. Precisely, $D_{k+1/2}(\tau)$ is the unique
polynomial of degree $4q+1$ satisfying for $\ell=0, \ldots ,2q$:
\be
\f {d^{\ell} D_{k+1/2}}{dt^{\ell}} (\tau_k) = \f {d^{\ell} P_{k}}{dt^{\ell}} (\tau_k) , \ \    
\f {d^{\ell} D_{k+1/2}}{dt^{\ell}} (\tau_{k+1}) = \f {d^{\ell} P_{k+1}}{dt^{\ell}} (\tau_{k+1}) . \label{Hermint}
\ee
As the interpolation operators are linear, and $D_{k+1/2}(\tau)$ depends only on the
stencil data, there exist degree-$(4q+1)$ polynomials  $\omega_{k+1/2,r}(\tau)$, $\tau_r \in S_k \cup S_{k+1}$,
such that
\be
D_{k+1/2}(\tau) = \sum_{r:\ \tau_r \in S_k \cup S_{k+1}} \!\! y_r \omega_{k+1/2,r}(\tau),
\qquad \tau\in (\tau_k,\tau_{k+1}).
\label{localexpD}
\ee
\item[iv.]
Let $k(\tau)$ be the index of the rightmost node not larger than $\tau$, i.e.\ such that
$\tau_{k(\tau)} \le \tau < \tau_{k(\tau)+1}$.
Then the D-spline basis function $\omega_r(\tau)$ associated with node $r$ is defined as the piecewise degree-$(4q+1)$ polynomial
\be
\omega_r (\tau) = \left\{ \ba{ll} \omega_{k(\tau)+1/2,r}(\tau), & \mbox{if $r$ is such that } \tau_r \in S_{k(\tau)} \cup S_{k(\tau)+1}, \\ 0, &
\mbox{otherwise} . \ea \right.
\label{Dsplinedef}
\ee
\end{description}

By construction $\omega_r(\tau) \in C^{2q}$, and, for uniformly spaced nodes in time,
they are translates of a
single simple basis function away from $\tau=0$. In addition, the D-spline interpolant
reproduces polynomials of degree $2q$, which by standard results implies accuracy of
order $2q+1$ for function values and $2q$ for derivatives.\footnote{Here we will identify the interpolants by the polynomial degrees they exactly reproduce; that is the degree-$2q$ D-spline is the piecewise degree $4q+1$ function described above.} 

\begin{figure}[htb]
\begin{center}
\hbox{
\includegraphics[width=.49\textwidth]{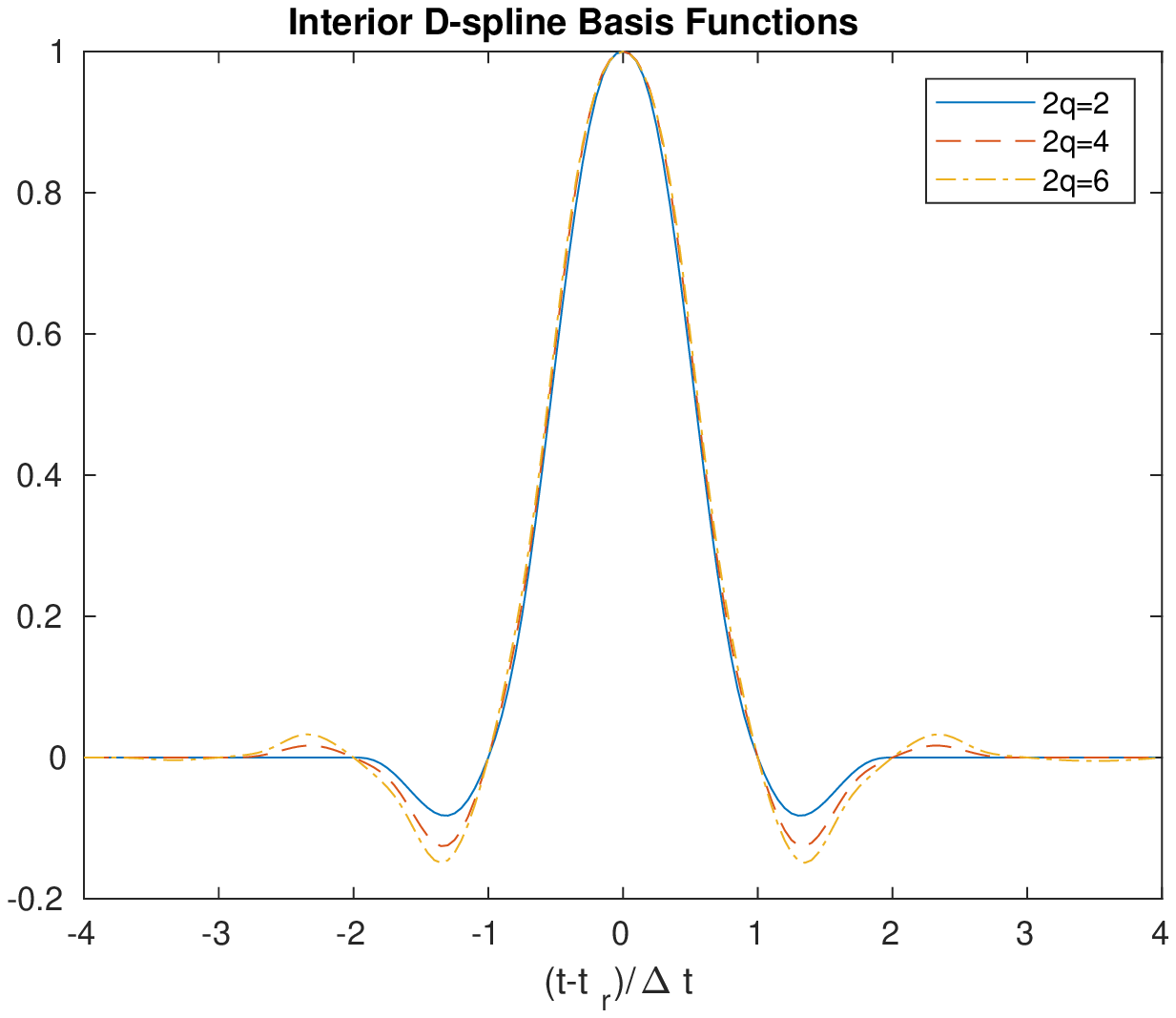}
\includegraphics[width=.49\textwidth]{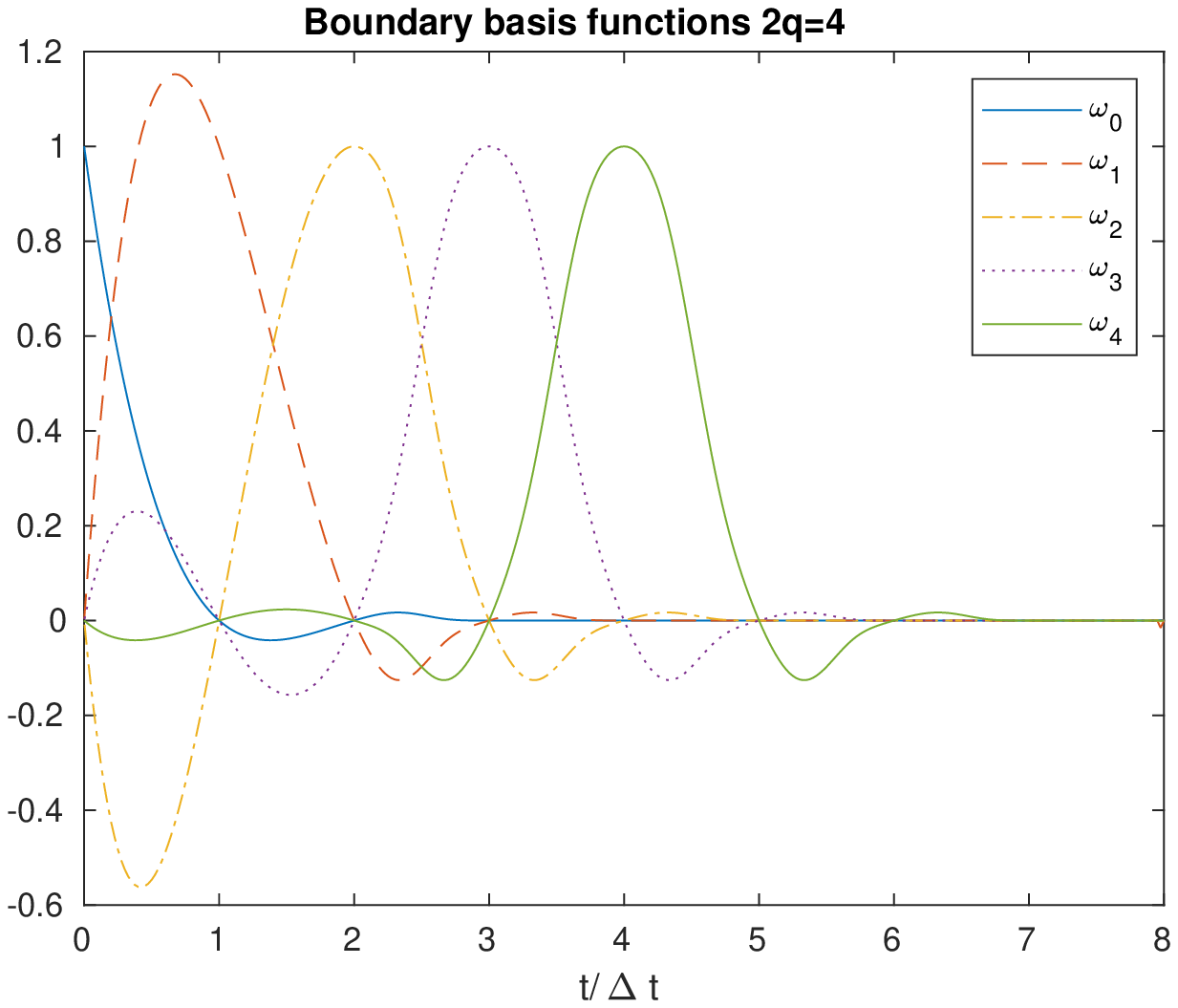} 
}
\caption{D-spline basis functions for the uniform grid case. (a) interior case, showing three different orders; (b) boundary case, showing basis functions for various nodes, at a fixed order.
\label{DsplineFig}}
\end{center}
\end{figure}

In the uniform grid case, we plot the
interior functions $\omega_r$ for
$2q=2$ to $6$, along with boundary functions for $2q=4$, in Figure \ref{DsplineFig}.  
A Fortran90 implementation (with MATLAB interface) of the above D-splines on regular grids is available
in the {\tt timedomainwaveeqn/timeinterp} directory; see Remark~\ref{r:code}.

\begin{figure}
\centering
\includegraphics[width=.49\textwidth]{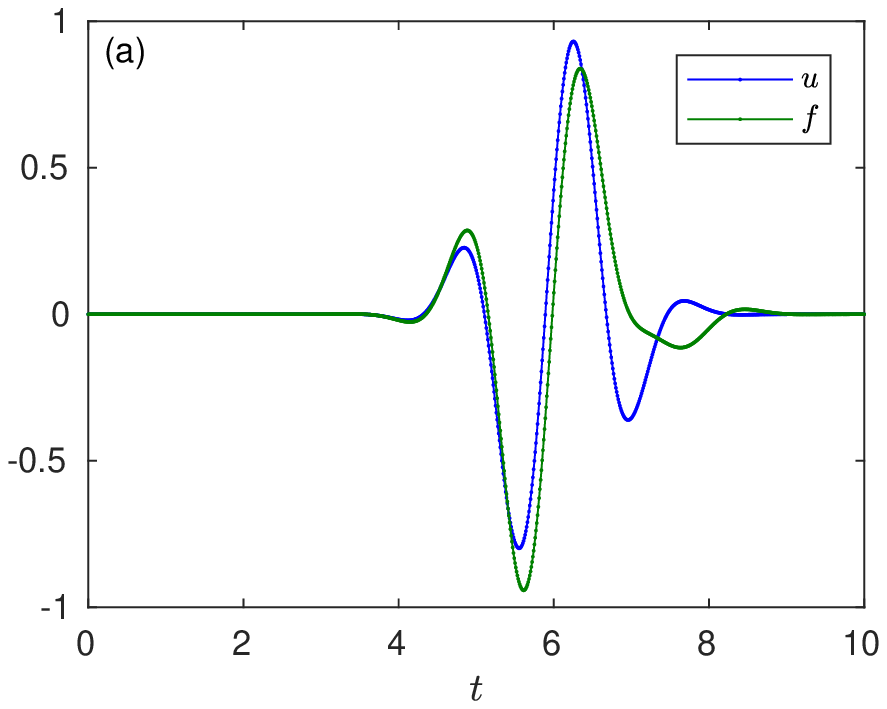}
\includegraphics[width=.49\textwidth]{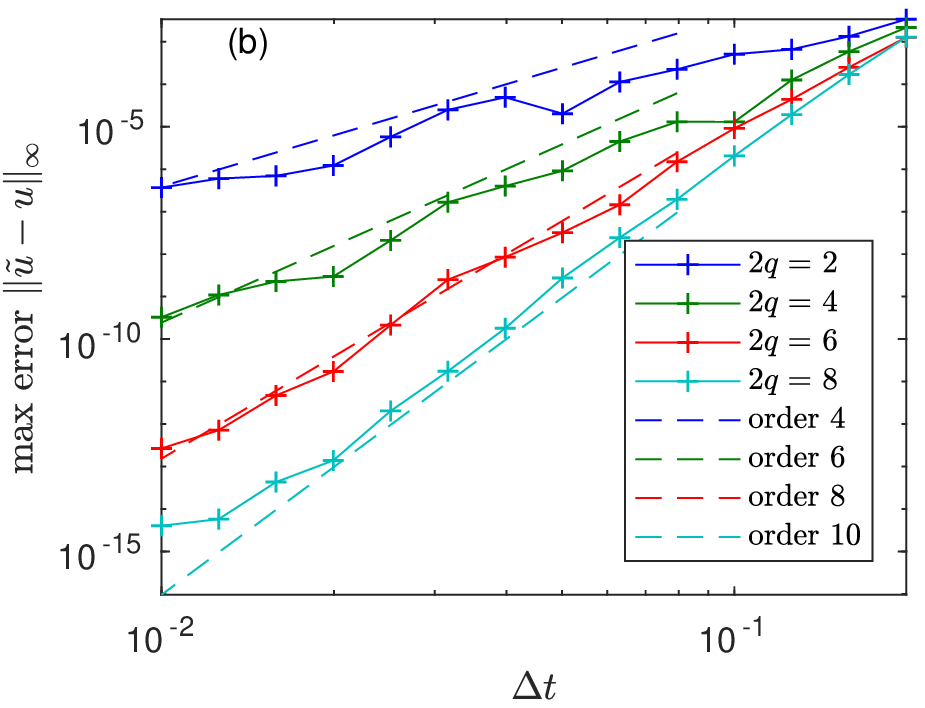}
\caption{Numerical solution of the Volterra equation \eqref{volt} over $t\in[0,10]$, with $f$ constructed so that the solution is $u(t)=e^{-(t-6)^2}\cos(4t)$. (a) Graph of $u(t)$ and right-hand side $f(t)$; the dots are spaced at the smallest tested $\Delta t = 0.01$. (b) Convergence with respect to timestep $\Delta t$ of the max error in $u$, for various spline basis orders $2q$, compared to convergence orders $2q+2$. We fix $p=16$.}
\label{f:volterra}
\end{figure}

\subsection{Application to a simple Volterra equation}
\label{s:simple}

Before attacking the full spatio-temporal problem,
we illustrate the use of the above temporal basis functions to create a
collocation
time-stepping scheme
for a simple 2nd-kind Volterra integral equation with ``top-hat'' kernel,
\begin{equation}
   u(t) + \int_0^1 u(t-s) ds \; = \; f(t)~, \qquad t\in[0,T]~.
\label{volt}
\end{equation}
We assume $f$ is smooth, and $u(t)$ is known for $t\le0$ and smooth for all $t\le T$.
As in Section~\ref{s:intro} let $t_k=k \Delta t$, $k\in\mathbb{Z}$
be the time grid.
Let $\{\tau_\ell\}_{\ell=1}^p$ be the nodes and $\{w_\ell\}_{\ell=1}^p$ be the corresponding weights of a $p$-node Gauss--Legendre quadrature rule on $[0,1]$.
Enforcing \eqref{volt} at $t=t_k$, and inserting the quadrature,
\begin{equation}
    u(t_k) + \sum_{\ell=1}^p w_\ell u(t_k-\tau_\ell) \;\approx\; f(t_k)~, \qquad k=1,\ldots,T/\dt~,
\label{colloc}
\end{equation}
holds to high accuracy because the integrand is smooth.
Let the unknowns be $u^k := u(t_k)$, for $k=1,\ldots,T/\dt$.
Let $\omega_r(\tau)$ be the degree-$2q$ D-spline basis, as defined in the
previous section, for the regular grid of spacing $\dt$,
which one may think of as stepping backwards in time from the current time $t_k$.
The resulting interpolant is
\begin{equation}
  u(t_k-\tau) \approx \sum_{r=0}^\infty \omega_r(\tau) u^{k-r}~,
  \label{interp}
\end{equation}
where, abusing notation slightly, we take $u^k$ for the pre-history $k\le 0$
to be populated with the known solution $u(t_k)$.
Substituting this interpolant into \eqref{colloc} defines a lower-triangular
Toeplitz linear system
\begin{equation}
    u^k + \sum_{r=0}^\infty W^r u^{k-r} = f(t_k)~, \quad k=1,\ldots,T/\dt~,
\label{toep}
\end{equation}
with weights computed by
\be
W^r = \sum_{\ell=1}^p w_\ell \omega_r(\tau_\ell)~,  \qquad r=0,1,\dots
\label{Wr}
\ee
Note that in \eqref{toep} the upper limit for $r$ can be replaced by $n\ge 1/\dt + q$, which is sufficiently large to capture all history dependence in
\eqref{volt}.

The system \eqref{toep} is naturally best solved sequentially for unknowns $k=1,2,\dots$, i.e.\ by time-stepping, since each row can be written
\begin{equation}
    (1 + W^0) u^k  \;=\;  f(t_k) - \sum_{r=1}^\infty W^r u^{k-r}  =: \tilde f^k~, \quad k=1,\ldots,T/\dt~,
\label{toepsol}
\end{equation}
where the ``right-hand side''
$\tilde f^k$ is explicitly given in terms of data and previous values.
Of course, in this simple example, the ``solve'' for $u^k$ at each time step
is trivial:
\begin{equation}
u^k = (1+W^0)^{-1} \tilde f^k~.
\label{impl}
\end{equation}

Figure~\ref{f:volterra} shows that this scheme achieves an empirical order $2q+2$.
In this experiment, $f$ is constructed numerically (via an accurate quadrature) such that the solution $u$ is a known function.
Note that, although in this test case $u(t)$ is not strictly zero for $t\le 0$,
it is of order $10^{-16}$ or less, so that initialization by zero data for the pre-history of $u$ induces negligible error.

The figure also shows stability for arbitrarily small $\dt$,
even with fixed number $p$ of quadrature nodes.
The mean mesh spacing $1/p$ is the closest analog to $\dx$ in the
full spatio-temporal scheme of Section~\ref{SpaceD}.
Thus, in contrast to the full scheme, and also to the modal problem for the
sphere in Section~\ref{modal}, we observe no inverse CFL condition.
In fact, for most of the $\dt$ tested, the quadrature \eqref{Wr} is completely unresolved for each basis function $\omega_r$, yet the scheme is accurate---since $u$ itself is resolved---and stable.

\begin{rmk}
Recall that in the context where there are $N$ spatial variables, the
analog of the solution \eqref{impl} to \eqref{toepsol}
requires the solution of an $N\times N$ linear system at each time step $k$,
i.e.\ an ``implicit'' time step, which can incur a large cost.
As mentioned above, for the full problem we will focus on 
predictor-corrector schemes, which replace 
this by a {\em fixed number} of explicit steps and require only $\bigO(N)$ effort.
\end{rmk}


\subsection{Modal problem on the unit sphere}
\label{modal}

As a second example of a simple Volterra-like scalar problem, but one which is more directly related to the full problem we aim to solve, we consider our combined field equation with $a=b=1$ restricted to the amplitudes, $\mu_n (t)$, of a spherical harmonic expansion of the full solution on the unit sphere. As shown in \ci{TDIEstab}, the integral equation to be solved is
\be
\f {\mu_n(t)}{2} + \f {1}{4} \int_0^2 P_n (1-s^2/2) \left( \mu_n (t-s) + (2-s) \mu_n' (t-s) \right) ds \;=\; g_n (t) , \label{spheren}
\ee
where $P_n$ denotes the degree-$n$ Legendre polynomial. The key difference between this example and the one considered above is the presence of the time derivative of the density. We believe this term leads to the inverse-CFL constraint observed for the full model, and want to verify that it appears in this simpler case where no spatial integrations are required. We note that in \ci{TDIEstab} the time derivative term is removed via integration by parts. Such a procedure would remove the difficulty, but it is unclear how it could be accomplished for the full problem with a general scatterer. 

\begin{figure}
    \centering
    \includegraphics[width=0.7\textwidth]{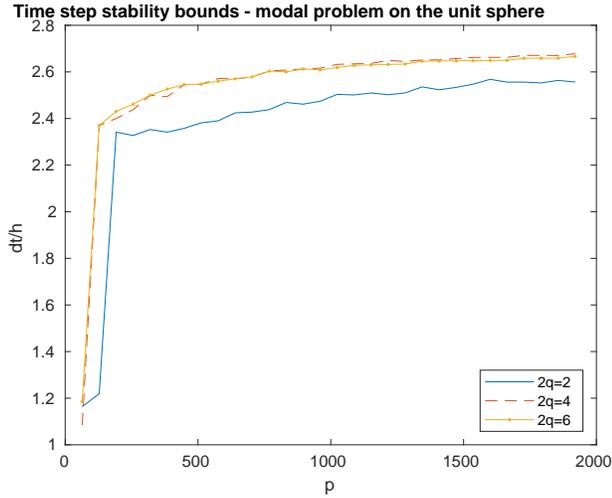}
    \caption{Minimum stable CFL number, $\dt/\dx$, versus number of quadrature nodes, $p$, for the modal problem on the unit sphere. Here $h=2/p$.}
    \label{Spherenfig}
\end{figure}

Our goal will simply be to determine the stability limits, if any, on the choice of time step. We discretize exactly as above: the integral is approximated by a $p$-point Gauss rule and the D-spline interpolant of $\mu_n$ as well as its derivative is evaluated at the quadrature nodes to produce a scheme
\be
\mu_n^k = (1/2+V^0)^{-1} \biggl( g_n (t_k) - \sum_{r=1}^{\infty} V^r \mu_n^{k-r} \biggr)~,
\ee
with weights
\be
V^r = \f {1}{4} \sum_{\ell=1}^p w_{\ell} P_n (1-\tau_{\ell}^2/2) 
\bigl( \omega_r (\tau_{\ell}) + 
(2- \tau_{\ell} )\omega_r' (\tau_{\ell}) \bigr)~.
\ee
Choosing $g_n$ to be, for every $n$ that we test,
\bd
g_n (t) = 10e^{-10 (t-2)^2} ,
\ed
we solve to $T=100$, for $p$ varying between $64$ and $1920$ in increments of $64$ and for all $n$ from $0$ to $p/2$. The method is deemed unstable for a given time step if the maximum value of the solution {\em for any $n$} grows beyond the expected value by a factor of roughly $1.1$. If we take the average mesh width to be $h = 2/p$ and plot $p$ versus the minimum stable CFL number, $\dt_{\tbox{min}}/h$, for $2q=2,4,6$ we obtain the results shown in Figure \ref{Spherenfig}. We see for $h$ small the minimum CFL number is approximately $2.6$ for $2q=4,6$ and $2.5$ for $2q=2$. We note that for $2q=8$ (a $10$th-order method) the stability region shrinks considerably since a maximum limit of about $3$ also appears; that is, the method is stable in this test only for $2.5h < \dt <3h$.
Thus we restrict ourselves to methods with $2q$ in the range 2 to 6
(orders $4$ to $8$), where the time step, once above its minimum stable value, appears to be controlled by accuracy and not stability, although we have no formal proof of such a stability result.

\begin{figure}
  \includegraphics[width=0.33\textwidth]{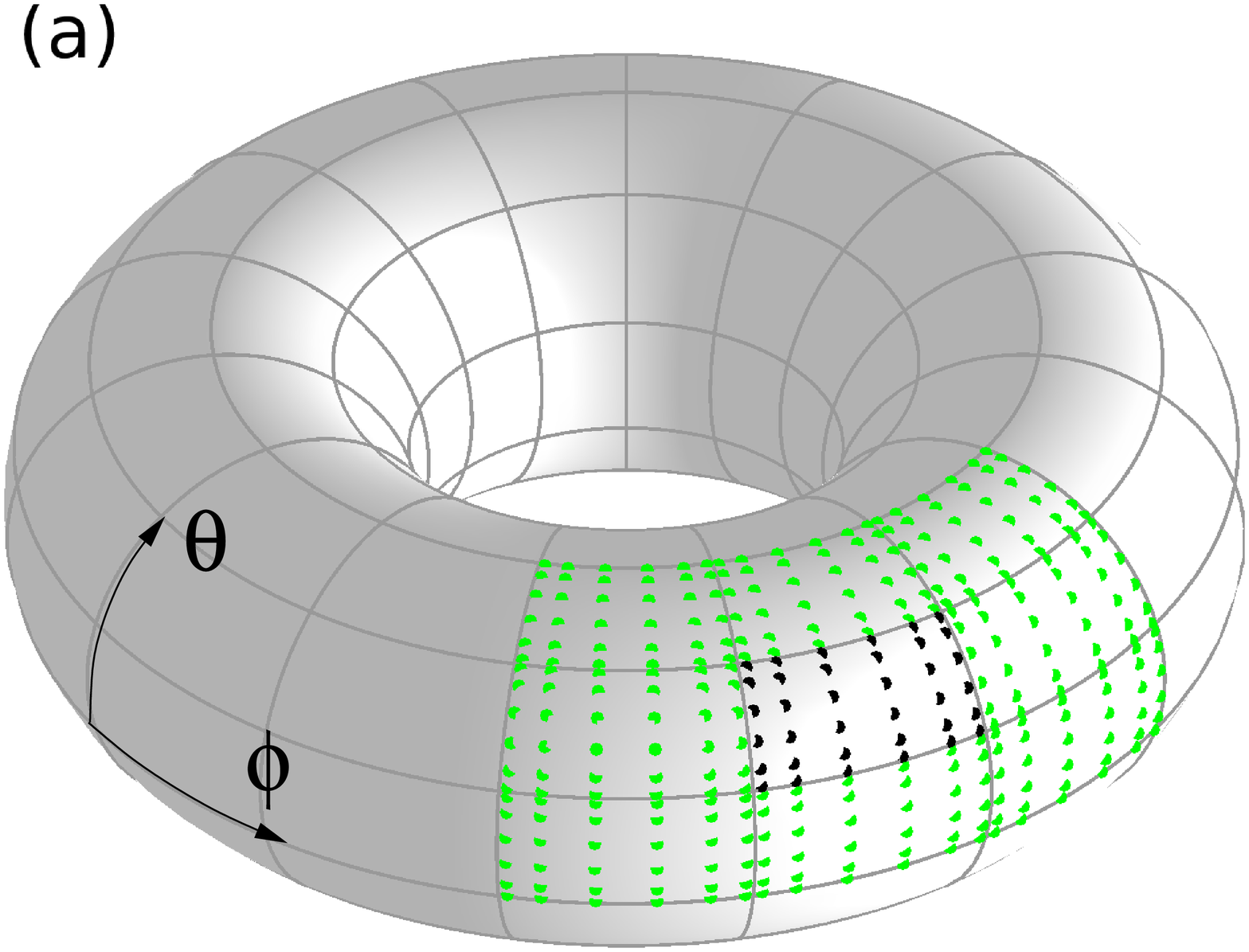}
  \includegraphics[width=0.33\textwidth]{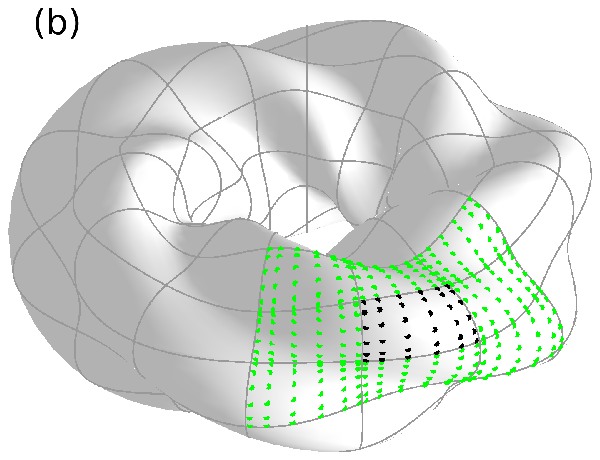}
  \quad
  \includegraphics[width=0.27\textwidth]{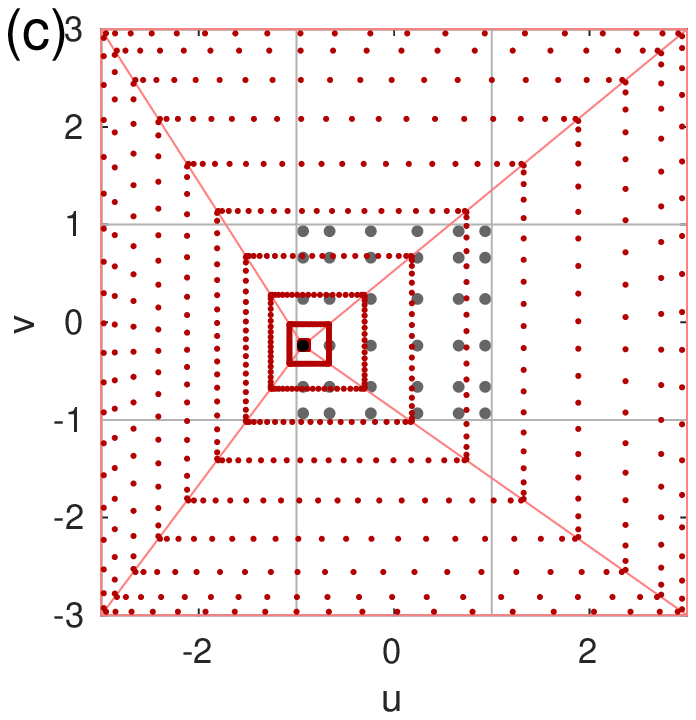}
  \caption{Surface $\Gamma$ and spatial singular quadrature scheme.
    (a) A torus surface parametrized by $(\phi,\theta)\in[0,2\pi)^2$,
      with $\nph=15$ by $\nth=10$ panels,
    showing panel divisions (grey lines), Nystr\"om nodes on one panel (black)
    and on its eight neighbors (green).
    ``Far'' nodes for the black nodes are not shown.
    (b) The ``cruller'' surface, showing the same; see section~\ref{NumExp}.
    (c) A standard panel $(u,v) \in [-1,1]^2$, with its eight neighbors, forming a
    $3\times 3$ grid (grey lines).
    Node preimages (grey) are shown only for the central panel.
    The auxiliary quadrature nodes (red) are shown with $n_r=n_\varphi/2=10$,
    for one target node $x_i$
    (preimage shown in black).
    These auxiliary nodes are polar Gauss--Legendre nodes over four triangles
    (outlined in red).%
    \label{f:geom}}
\end{figure}

\section{Spatial discretization and full scheme}\label{SpaceD}

First we present then test
a quadrature scheme to apply the retarded single- and
double-layer potentials \eqref{Slp} and \eqref{Dlp}.
Finally we present the full interpolation scheme and
Volterra time-stepping for \eqref{cftdie}.

\subsection{Retarded layer potentials for off-surface targets}
\label{s:offsurf}

The case of an exterior evaluation target $x\in\Omega$, far
from $\Gamma$, is simple:
for densities $\mu(y,t)$ that are smooth with respect to both $x\in\Gamma$
and time $t$, the integrands are also smooth,
and a standard quadrature scheme using the density interpolatory nodes
$x_j\in\Gamma$ will be accurate.
In this work we restrict ourselves to smooth torus-like surfaces.
Such surfaces can be parameterized by an infinitely differentiable,
doubly $2\pi$-periodic function $x=z(\phi,\theta)$,
where $z:[0,2\pi)^2 \to \RR^3$.
We now describe a simple composite (i.e.\ panel-based) rule to generate
the nodes $x_j$ and weights $w_j$, which we use in later tests,
such that for any smooth $f:\Gamma\to\RR$,
\be
\int_\Gamma f(y) dS_y = \int_0^{2\pi}\!\!\int_0^{2\pi}\! f(\phi,\theta)\,
\|z_\phi\times z_\theta\|\,d\phi d\theta
\; \approx \; \sum_{j=1}^N w_j f(x_j)
\label{quad}
\ee
holds to high-order accuracy.
Here $z_\phi$ and $z_\theta$ are the partials of $z$.
We cover the parameter space $[0,2\pi)^2$ with a uniform
$\nph$-by-$\nth$ grid of rectangular patches (``quads'').
Each patch is covered by a tensor product grid comprising a
$p$-node Gauss--Legendre rule in each of the two directions.
Let $x_j$, $j=1,\dots,N$, where $N=\nph\nth p^2$,
be the images of these parameter nodes under the map $z$.
Figs.~\ref{f:geom}(a) and (b) show two surfaces with some of their resulting
nodes.
The corresponding weights $w_j$ are found as follows.
Let $\eta_m$, $m=1,\dots,p$, be Gauss--Legendre weights
on the interval $[0,1]$.
For surface node $j$, let $\phi_j$ and $\theta_j$ be its parameter values, and
$m_j, m'_j \in \{1,p\}$ be its indices in the two directions within the
appropriate panel.
Then
\be
w_j = \frac{(2\pi)^2}{\nph\nth} \eta_{m_j}\eta_{m'_j}\,
\| z_\phi(\phi_j,\theta_j) \times z_\theta(\phi_j,\theta_j) \|~.
\label{wj}
\ee
The expected convergence
order for \eqref{quad} is $\bigO(h^{2p})$,
where $h=\bigO(N^{-1/2})$ is the resolution
(combining \cite[(2.7.12)]{davisrabin} with a theorem on composite rules
\cite[Sec.~2.4]{davisrabin}).

\subsection{Retarded layer potentials for on-surface targets}
\label{s:onsurf}

When the target $x$ is on $\Gamma$, as needed
in the integral equation \eqref{cftdie},
then the integrand has the following type of weak singularity.
If one smoothly parametrizes $\Gamma$ via local polar coordinates
$(r,\varphi)$ centered at the target point $x$,
for each $\varphi$ (i.e.\ radial line)
the integrand is $1/r$ times a smooth function of $r$.
Surprisingly, this is
the same form as for the 3D elliptic BVP case:
the kernels in \eqref{Slp} and \eqref{Dlp} are identical to
(or in the 2nd term of \eqref{Dlp}, less singular than) the Laplace
kernels, and although the retardation introduces a conical singularity
into the density, this does not change the singularity of their product.

For the elliptic case there exist many high-order Nystr\"om
quadrature schemes in 3D.
For surfaces diffeomorphic to the sphere,
a global
spherical harmonic basis \cite{wienert90} \cite[Sec.~3.6]{coltonkress},
or spherical grid rotation
\cite{ganesh,gimbutasgrid}
achieves spectral accuracy.
For more general
smooth surfaces, 
Bruno--Kunyansky \cite{brunoFMM,ying06}
use a smooth partition of unity to isolate the singular
near-target contribution.
To handle the latter they exploit the fact that
the polar metric $r dr d\varphi$ cancels the $1/r$ singularity
in the integrand, so integrate using auxiliary quadrature nodes on a polar grid
centered at the target.
For general high-order triangulations, including those with high aspect ratio,
Bremer--Gimbutas \cite{bremer3d}
developed generalized Gaussian quadratures
that again use auxiliary nodes.
Note that in these elliptic schemes,
the integral kernel is directly evaluated at each of the auxiliary nodes,
but the density must be spatially interpolated
from its values at the nodes $x_j$.

Building on the above, we present a simple high-order accurate scheme to
evaluate the retarded layer potentials \eqref{Slp} and \eqref{Dlp}
on surfaces discretized
by a structured rectangular grid of quad panels of the type described
above in Sec.~\ref{s:offsurf}.
Let $x_i$ be a target node, in panel $k$.
This panel and its eight neighbors form a $3\times 3$ block of ``near'' panels,
containing node indices $j\in \Jni$ (the
black and green nodes in Fig.~\ref{f:geom}(a)),
leaving $(\nph\nth-9)$ ``far'' panels.
With respect to each of the latter panels, the singularity of the kernel
is distant, so their native quadrature nodes $x_j$ and weights \eqref{wj}
may be accurately used,
as in the previous section.
This explains the first term in our approximation
\be
\int_{\Gamma} \frac{f(y)}{4\pi |x_i-y|}dS_y
\;\approx\;
\sum_{j\notin\Jni} \frac{w_j}{4\pi r_{ij}} f(x_j)
+
\sum_{\ell=1}^{N_\tbox{aux}} \tilde{S}_{i\ell} f(y_{i\ell})
~,
\label{Ssplit}
\ee
where $r_{ij} := |x_i-x_j|$;
the second term is explained below.
Note that \eqref{Ssplit} is of the form \eqref{matrices},
with $N'=(\nph\nth-9)p^2+N_\tbox{aux}$.
The expressions for the other two kernels are analogous.
  
The second term in \eqref{Ssplit}
accounts for the contribution of the near $3\times 3$ panel block
via a new target-specific set of $N_\tbox{aux}$ auxiliary nodes.
For this it is convenient to switch to
the standard parametrization $(u,v)\in[-3,3]^2$
for this near block, with the target panel $k$ preimage being $[-1,1]^2$;
see Fig.~\ref{f:geom}(c).
Precisely, there is a simple affine map to the global parameters
$\phi = \pi(2i_k +u+1)/\nph$ and
$\theta = \pi(2i'_k +v+1)/\nth$,
where $i_k$ and $i'_k$ are the
toroidal and poloidal indices (i.e.\ integer coordinates) of panel $k$
within the panel grid.
Then denote by $\tilde z$ the ($k$-dependent) map from $(u,v)$ to $\RR^3$, which
is the above affine map composed with the map $z$.
The auxiliary nodes comprise four grids, each of which
integrates over one of the four triangles connecting
the block walls to the preimage of the target.
Together the four triangles cover $[-3,3]^2$; see Fig.~\ref{f:geom}(c).

Specifically,
let the polar coordinates $(r,\varphi)$ be centered at the target preimage
in $(u,v)$,
and consider one of the four triangles, $T$, that
lies in the angle range $[\varphi_0,\varphi_1]$ and whose far edge
is given by $r(\varphi)$ in polar coordinates.
For the single-layer \eqref{Slp},
the integral of a retarded density $f$ as in \eqref{f} over the
image of this triangle on $\Gamma$ is
\be
\iint_T \frac{J(u,v) f(u,v)}{4\pi |x_i - \tilde z(u,v)|}
du dv = \int_{\varphi_0}^{\varphi_1}
\int_{0}^{r(\varphi)} \frac{J(r,\varphi) f(r,\varphi)}
{4\pi |x_i - \tilde z(r,\varphi)|}
\, rdr \, d\varphi
\label{tri}
\ee
where $J(u,v) = \| \tilde z_u(u,v) \times \tilde z_v(u,v)\|$
is the Jacobian
of the map to the surface,
and we abuse notation slightly so that $J(r,\varphi)$ means
$J(u(r,\varphi),v(r,\varphi))$, etc.
Let $\rho_m$ and $\eta_m$, $m=1,\dots,n_r$, be respectively
the nodes and weights
of a $n_r$-point Gauss--Legendre rule on $[0,1]$.
Let $\varphi_n$ and $\xi_n$, $n=1,\dots,n_\varphi$, be the nodes and weights of a $n_\varphi$-point Gauss--Legendre rule on $[\varphi_0,\varphi_1]$.
Then \eqref{tri} is approximated by
\be
\sum_{n=1}^{n_\varphi} \xi_n [r(\varphi_n)]^2 \sum_{m=1}^{n_r}
\eta_m \frac{J_{nm} f_{nm}}{4\pi|x_i - \tilde z_{nm}|} \rho_m~,
\label{auxrule}
\ee
where $J_{mn} := J(\rho_m r(\varphi_n),\varphi_n)$, etc, indicates
the value at the auxiliary node indexed by $n$ and $m$.
High-order convergence is expected for \eqref{auxrule} since,
although $f$ has a conical singularity at the polar origin,
along constant-$\varphi$ rays the integrand times $r$ is smooth.
Summing the four triangles,
there are $N_\tbox{aux} =4n_rn_\varphi$ auxiliary nodes for each
target point.
The weights $\tilde S_{i\ell}$ in \eqref{Ssplit} may be read off
by associating each $\ell$ with a term $nm$ in \eqref{auxrule},
and taking all factors except the density sample $f_{nm} = f(y_{i\ell})$.
The weights $D_{i\ell}^h$ and $W_{i\ell}^h$ in \eqref{matrices} are
found in an analogous way.

\begin{rmk}[order of convergence]
  The convergence of this on-surface scheme is subtle,
  due to its split into far and near source panels.
  Consider the error due to the $p\times p$-node smooth rule on one of
  the {\em nearest} ``far'' panels (i.e.\ just outside of the $3\times 3$
  block):
  as $h\to 0$, its integrand has a singularity that remains at a {\em fixed
  distance relative to the panel size}, so its error is expected to drop
  only in proportion to the panel area times the typical integrand.
  Thus, for any fixed $p$, the formal order $\gamma$ is
  low ($\gamma<2$).
  Yet, this matters little in practice 
  because its
  prefactor is expected to be exponentially small in $p$.
  This follows by analogy with 1D $p$-node Gauss quadrature on $[-1,1]$,
  for which the error is of order $\varepsilon = \rho^{-2p}$,
  where $\rho$ is the size parameter of a Bernstein ellipse in which the integrand
  is analytic \cite[Thm.~19.3]{ATAP}.
  Since the nearest singularity is at $\pm 3$ scaled to the standard
  panel $[-1,1]$, solving $3 = (\rho+\rho^{-1})/2$ gives
  $\rho \approx 5.8$, hence for $p=4$, $\varepsilon < 10^{-6}$,
  and for $p=8$, $\varepsilon < 10^{-12}$.
  Combining with the accuracy of the rest of the far panels, and assuming
  $n_r > p$,
  one can summarize the expected error as
  $\bigO(\varepsilon h^q + h^{2p})$.
  We postpone a rigorous analysis for future work.
  \label{r:order}
\end{rmk}

\begin{rmk}
  Although similar to that of \cite{bremer3d}, our method is simpler
  and somewhat more efficient,
since we exploit the structured nature of the panel grid
to cover both self-interaction and neighboring panel interactions
with a single auxiliary node rule.
This is only possible because, for our
class of surfaces, the charts from each panel
extend to their neighbors in a known, smooth fashion.
We also do not attempt to handle panels of aspect ratios much larger than 2.
\end{rmk}

\begin{figure}
  \mbox{\hspace{-2ex}
    \raisebox{0.7in}{
      \includegraphics[width=0.26\textwidth]{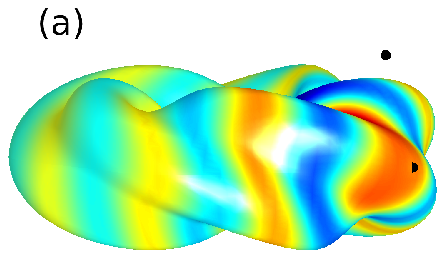}}\;
  \includegraphics[width=0.4\textwidth]{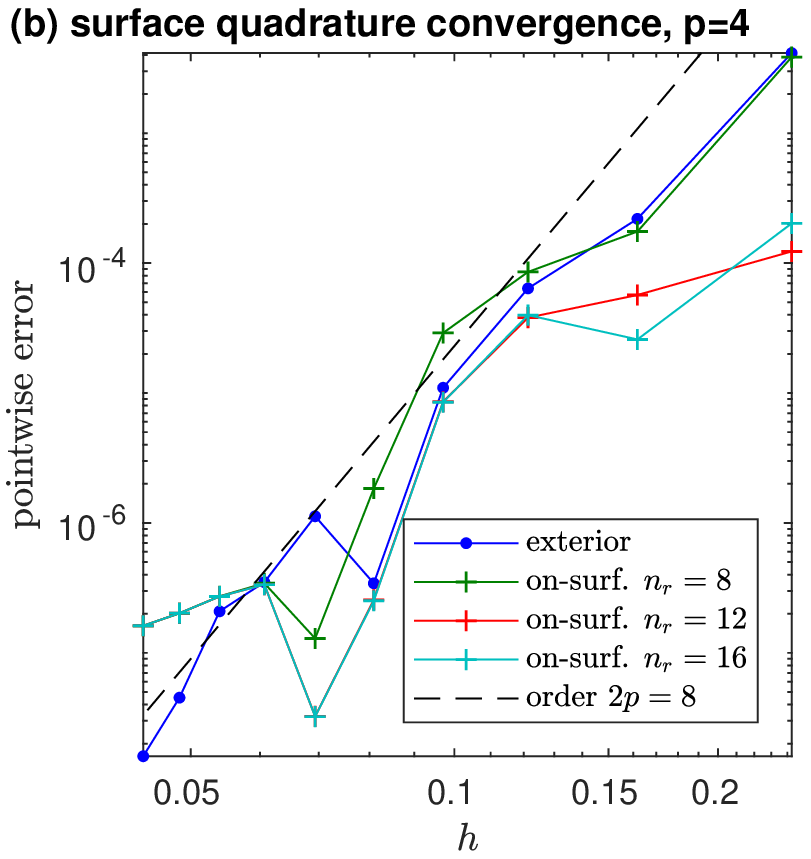}
  \hspace{-4ex}
  \includegraphics[width=0.4\textwidth]{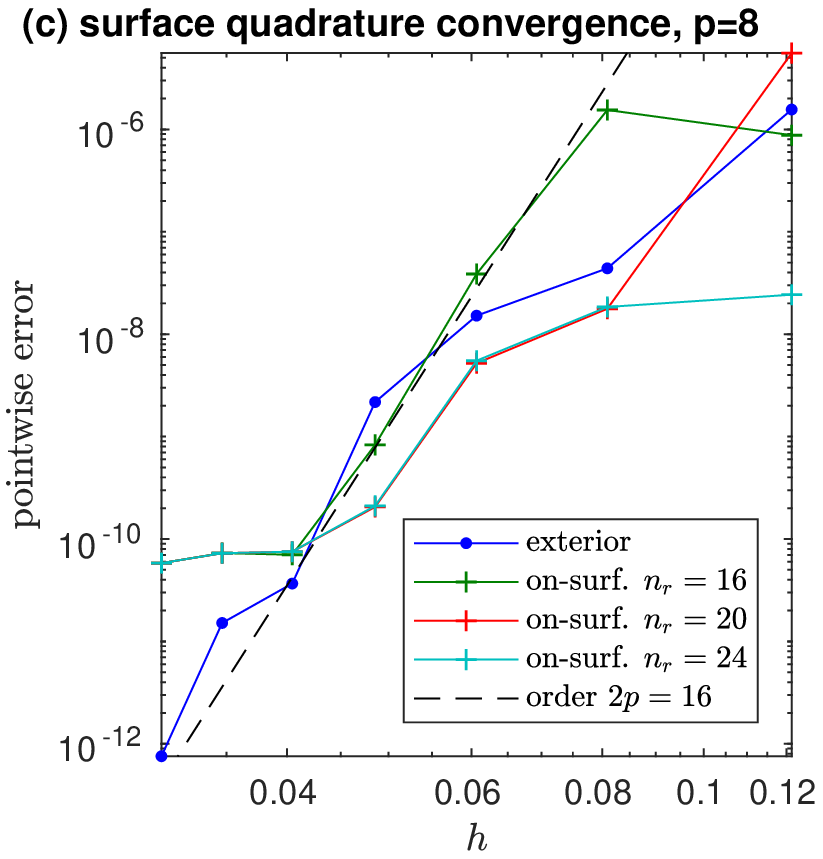}}
  \vspace{-4ex}
  \caption{Test of surface quadrature scheme via Green's representation
    formula for the cruller of Fig.~\ref{f:geom}(b).
    (a) The dots show the off-surface and on-surface test targets $x$. The
    colors show for $y\in\Gamma$
    the retarded double-layer density $u^+(y,t-|x-y|)$ for $x$
    on-surface; note this has a conical singularity at $x$.
    (b) Errors for an exterior target and on-surface target,
    for panels of order $p=4$.
    (c) Same for order $p=8$.
    In the last two panels, $\dx$ is defined by \eqref{dxdef}, and
    various auxiliary node orders
    $n_r=n_\varphi/2$ are compared for the on-surface case;
    see Sec.~\ref{s:onsurf}.
  \label{f:spaceconv}}
\end{figure}  

\subsection{Validation of retarded layer potential evaluation}
\label{s:valid}

Before presenting the full spatio-temporal scheme, we pause to
describe some numerical tests of the above spatial quadratures
for retarded potentials.

{\bf Surfaces used for tests.}
In this work we use a simple class of smooth surfaces $\Gamma$ diffeomorphic
to a torus.
Given $H(\phi,\theta)$, a doubly $2\pi$-periodic smooth height modulation
function,
then the map $x=z(\phi,\theta)$ is given in Cartesian coordinates by
\be
z(\phi,\theta) = \bigl(
(1 + H(\phi,\theta)\cos \theta)\cos \phi, \,
(1 + H(\phi,\theta)\cos \theta)\sin \phi, \,
H(\phi,\theta) \sin \theta \bigr)~.
\label{torus}
\ee
Thus, the major radius is 1.
The plain torus of Fig.~\ref{f:geom}(a) is given by the constant $H\equiv 0.5$,
and is close to being the most benign surface of this topology.
The ``cruller'' of Fig.~\ref{f:geom}(b)
has $H(\phi,\theta) = 0.5 + 0.1 \cos(5\phi + 3\theta)$,
and is more challenging due to its higher curvature,
and ridges which do not align with the parameter directions.
We will not list the analytic partial derivatives of $z$ here,
although of course they are needed for Jacobian computations.
Note that both objects have a maximum diameter of approximately 3.

{\bf Green's representation formula.}
We test the convergence of the
above quadrature schemes by numerically checking Green's representation formula
(Kirchhoff's formula) for the wave equation
\cite[(1.17)]{TDIE}, which for convenience we now state.
Let $u(x,t)$ satisfy \eqref{waveq} in the closure of the
exterior domain $\Omega$, for all $t$,
and let $u^+$ and $u_n^+$ indicate respectively the value and
normal derivative on $\Gamma$,
then
\be
({\mathcal D} u^+)(x,t) - ({\mathcal S} u_n^+)(x,t) =
\left\{
\begin{array}{cc}
  u(x,t), & x \in \Omega,\\
  u(x,t)/2, & x \in \Gamma.
  \end{array}
\right.
\label{GRF}
\ee
(Note that, since $u$ is a solution for all time $t$,
there is no need for initial conditions as in \cite[Sec.~1.4]{TDIE}.)
We use this to test both the off-surface (Sec.~\ref{s:offsurf})
and on-surface (Sec.~\ref{s:onsurf}) quadrature schemes.
We use for $u$ an exterior solution  to the wave equation generated by a
generic point source
$x_0 = (0.9,-0.2,0.1)$ inside, but far from, the surface $\Gamma$,
namely
\be
u(x,t) = T(t-|x-x_0|), \qquad x\in\Omega, \; t\in\RR,
\label{usrc}
\ee
with a source signal $T(t) = \cos 5t$.
Fig.~\ref{f:spaceconv}(a) shows a snapshot of the resulting retarded double-layer density
for an on-surface target $x$.
In this section we show results only for the cruller, omitting
the more accurate and predictable
results for the plain torus.

{\bf Off-surface test.}
Figs.~\ref{f:spaceconv}(b) and (c) shows (with dots) the convergence
of \eqref{GRF} for an exterior target $x = (1.3,0.1,0.8)$,
which is a generic point a distance 0.36 from $\Gamma$ (see subfigure (a)).
We use the smooth panel scheme of Sec.~\ref{s:offsurf}.
For both shapes, $\nth = 3\nph/2$ gives panels with low
aspect ratios. 
Recalling that each quadrature panel has $p\times p$ nodes,
(b) shows $p=4$ (for $N$ ranging from 384 to 11616), 
while (c) shows $p=8$ (for $N$ from 1586 to 24567).  
The horizontal axis shows the mean spatial node spacing, or resolution,
which we define as
\be
\dx := \biggl(\frac{\mbox{area}(\Gamma)}{N} \biggr)^{1/2}~.
\label{dxdef}
\ee
Although there is variation,
convergence is asymptotically consistent with the expected order $2p$.
The variation is absent and the errors much smaller for the plain torus
(not shown).
For $p=8$, 12-digit accuracy is reached at the highest $N$.

{\bf On-surface test.}
A surface target $x\in\Gamma$ is shown in Fig.~\ref{f:spaceconv}(a),
and the singular auxiliary node scheme of Sec.~\ref{s:onsurf}
used.
This target was in the corner of a panel, close to parameters $\phi=\theta=0$,
and thus involved worst-case auxiliary triangle aspect ratios.
For panel aspect ratios up to around 2, we found for
the auxiliary quadrature
that $n_\varphi=2n_r$ was adequate, so we fixed this ratio.
We explore three choices of $n_r$ for each choice of $p$.
For $p=4$, Fig.~\ref{f:spaceconv}(b) shows that,
for $h\le0.12$ ($\nph\ge 12$),
errors due to the singular scheme are negligible
relative to the overall error for $n_r=12$.
For $p=8$, 
a higher $n_r$ is needed: $n_r=20$ causes negligible errors for 
$h\le0.1$ ($\nph\ge 9$).
In general, for all but the largest $\dx$ available, the 
errors of the singular scheme are negligible for the choice $n_r=2p+4$.


\begin{rmk}[Choice of singular scheme order $n_r$]
  Based on such tests, we fix
  $n_r=2p$ for the torus and $n_r=2p+4$ for the cruller.
  Thus the auxiliary scheme is of much higher order than the underlying panels.
  This allows us to explore larger $h$ without loss of accuracy or
  stability. At smaller $h$ one could reduce $n_r$ and hence the effort
  for the auxiliary scheme.
  \label{r:nr}
\end{rmk}

With $n_r$ converged as above, the
$h$-convergence is as expected from Remark~\ref{r:order}:
errors drop with order roughly $2p$ (dashed lines), until they saturate
to a low-order convergence at around $10^{-7}$ (for $p=4$) or
$10^{-10}$ (for $p=8$).


\subsection{Interpolation and explicit time-stepping}
\label{s:predcorr}

We now describe how the above spatial quadrature for retarded potentials
is combined with the time-stepping of Sec.~\ref{TempD} to solve
the time-dependent BIE \eqref{cftdie}.
Enforcing the BIE on the time grid $t=t_k$ (as in Sec.~\ref{s:simple}),
and on the spatial nodes $x_i$, gives
\be
\frac{\mu^k_i}{2} +
\bigl[D\mu + S\bigl(a \frac{\pa\mu}{\pa t} + b\mu\bigr) \bigr](x_i,t_k)
= g^k_i,
\quad i=1,\dots,N, \quad k=1,\dots,T/\dt~.
\label{biecolloc}
\ee
Now, fixing constant $a$ and $b$,
we approximate the action of the retarded integral operators
on the density interpolated from the space-time data $\mu_j^k$
by a set of $N\times N$ matrices $A^r$, thus
\be
\frac{\mu^k_i}{2} + \sum_{r=0}^n \sum_{j=1}^N A^r_{ij} \mu^{k-r}_j
\;=\; g^k_i,
\qquad i=1,\dots,N, \quad k=1,\dots,T/\dt~.
\label{A}
\ee
where, analogously to Sec.~\ref{s:simple}, $n\ge \mbox{diam}(\Gamma)/\dt+q$
is large enough to capture all history dependence
via Huygens' principle plus the support of the D-splines.

{\bf Spatio-temporal interpolation.}
Each above matrix $A_r$ is filled as follows. For simplicity, 
consider only the plain single-layer ($S\mu$) contribution
in \eqref{biecolloc}.
Applying spatial quadrature \eqref{Ssplit}, then time interpolation \eqref{ConvApp}, gives
\bea
(S\mu)(x_i,t_k) &=&
\int_\Gamma \frac{\mu(x_i,t_k-|x_i-y|)}{4\pi|x_i-y|}
dS_y \nonumber \\
&\approx&
\sum_{j\notin\Jni} \frac{w_j}{4\pi r_{ij}} \mu(x_j,t_k-r_{ij}) +
\sum_{\ell=1}^{N_\tbox{aux}} \tilde S_{i\ell} \mu(y_{i\ell},t_k-|x_i-y_{i\ell}|)
\nonumber \\
&\approx&
\sum_{r=0}^n \sum_{j\notin\Jni} \frac{w_j}{4\pi r_{ij}} \omega_r(r_{ij})\mu_j^{k-r} +
\sum_{r=0}^n \sum_{\ell=1}^{N_\tbox{aux}} \tilde S_{i\ell} \omega_r(|x_i-y_{i\ell}|)\mu(y_{i\ell},t_{k-r})~.
\label{Sinterp}
\eea
The first term ($j$ far from $i$) is already in the form \eqref{A},
so the contribution to $A^r_{ij}$ in this case is
$\frac{w_j}{4\pi r_{ij}} \omega_r(r_{ij})$.
We now use spatial interpolation on each time slice $t_{k-r}$
to turn the auxiliary term also into a weighted sum over $\mu_j^{k-r}$.
Let $L_j(u,v)$, for $j\in\Jni$, 
be a set of $9p^2$ basis functions that interpolate
over the preimage of the near $3\times 3$ patch.
I.e., for any smooth function $f$,
\be
f(u,v) \;\approx \sum_{j\in\Jni} f_j L_j(u,v)~,
\qquad   (u,v)\in[-3,3]^2~,
\label{Lag}
\ee
holds to high order, where $f_j$ are the values at the interpolatory nodes
(preimages of $x_j$).
Let $(u_{i\ell},v_{i\ell})$ be the standard parameters of the $\ell$th auxiliary
node for the target node $i$.
Spatial interpolation of $\mu(\cdot,t_{k-r})$ then approximates the 2nd term of
\eqref{Sinterp} by
$$
\sum_{r=0}^n \sum_{\ell=1}^{N_\tbox{aux}} \tilde S_{i\ell} \omega_r(|x_i-y_{i\ell}|)
\sum_{j\in\Jni} L_j(u_{i\ell},v_{i\ell}) \mu_j^{k-r}~,
$$
which is now also of the form \eqref{A}.
Proceeding as above with the other two terms of \eqref{biecolloc}
(and recalling that $D$ has two terms \eqref{Dlp}), finally gives the formula
\be
A^r_{ij} = \left\{\begin{array}{ll} \frac{w_j}{4\pi r_{ij}}
\big[
  \bigl(\frac{n_j\cdot(x_i-x_j)}{r_{ij}^2} + b\bigr) \omega_r(r_{ij}) +
   \bigl( \frac{n_j\cdot(x_i-x_j)}{r_{ij}} + a\bigr) \omega'_r(r_{ij})
  \bigr],
& j\notin\Jni
\\
\sum_{\ell=1}^{N_\tbox{aux}}  L_j(u_{i\ell},v_{i\ell}) \bigl[
  (\tilde D_{i\ell} + b\tilde S_{i\ell} )\omega_r(|x_i-y_{i\ell}|) \\
  \hspace{1.5in}
  +\; (\tilde W_{i\ell} + a\tilde S_{i\ell} )\omega'_r(|x_i-y_{i\ell}|)
  \bigr],
& j\in\Jni
\end{array}\right.
\label{Arij}
\ee
For the basis $L_j$ we use the $p\times p$ product
Lagrange basis for whichever of the nine panels $j$ lies in, and zero
elsewhere.
Precisely,
let $k_j\in\{1,\dots,9\}$ be the panel in which node $j$ lies,
let $i_j$ and $i'_j$ be its two index coordinates within that panel,
and let $(u_0^k,v_0^k)$ be the parameter offset of panel $k$ relative
to the target panel ($u_0^k$ and $v_0^k$ are either $-2$, $0$, or $2$).
Then,
$$
L_j(u,v) = \left\{\begin{array}{ll} l_{i_j}(u-u_0^{k_j}) l_{i'_j}(v-v_0^{k_j})~,
& |u-u_0^{k_j}|\le 1, \,|v-v_0^{k_j}|\le 1~,\\
0~, & \mbox{otherwise~,}
\end{array}
\right.
$$
where $l_i(x)$ are the usual 1D Lagrange polynomials for the $p$
Legendre nodes on $[-1,1]$.
This form aids bookkeeping since it decouples all interactions into
independent panel-panel pairs, each of which is either near or far.
The $r$ different $p^2$-by-$p^2$ blocks of $A^r_{ij}$ given by all
nodes $j$ in a single source panel interacting with all nodes $i$ in
a single target panel may be filled together.
In practice the $N_\tbox{aux}\times p^2 \times p^2$ nonzero
entries of $L_j(u_{i\ell},v_{i\ell})$ are precomputed once and for all,
then for each target $j$ the sum over $\ell$ in \eqref{Arij} is performed as an
efficient matrix-matrix multiplication (GEMM).
Note that the near-panel interpolation error is expected to be $\bigO(\dx^p)$.

\begin{rmk}
  One might be tempted to interpolate with respect to $y$ the
  retarded densities such as $\mu(y,t-|x_i-y|)$; however,
  this would fail to be high-order accurate due to
  the conical singularity around the target $y=x_i$.
  Instead one must interpolate in {\em both} space and time, as above, since
  as a function of space and time $\mu(y,t)$ is smooth.
\end{rmk}

\begin{figure}  
  \includegraphics[width=\textwidth]{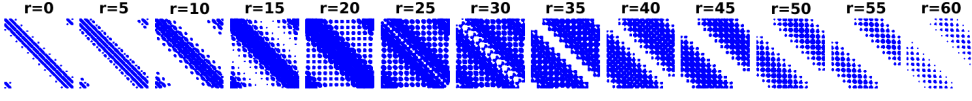}
  \\
  \vspace{-6ex}
\caption{Sparsity patterns of matrices $A^r$ for an arithmetic
  sequence of discrete time delay $r$ values. Each square shows
  $i$ down and $j$ across, with a node ordering fast within panels, slow
  between panels. Each matrix is $N\times N$ with $N=3456$.
  The average sparsity is $0.058$.
  The torus surface of Fig.~\ref{f:geom}(a)
  was used, with $\nph=12$, $\nth=8$, $p=6$, $2q=4$,
  and $\dt=0.05$.
}
\label{f:Arij}
\end{figure}

Fig.~\ref{f:Arij} shows the sparsity patterns of a selection of
the resulting $A^r$ matrices.
As expected by Huygens' principle,
$A^0$ is concentrated around the diagonal, but as the time delay
$r$ increases, the shell of influence spreads across the panels,
departing at the most distant (furthest off-diagonal) panels.


{\bf Predictor-corrector scheme.}
Using the notation $\bmu^k := \{\mu_j^k\}_{j=1}^N$ for the vector of
densities at time step $k$,
we rewrite \eqref{A} as a linear system to be solved at each time step,
\be
(\sfrac{1}{2}I + A^0)\bmu^k \; = \; \mbf{g}^k - \sum_{r=1}^n A^r \bmu^{k-r}
=: \tilde{\mbf{g}}^k
~, \qquad k=1,\ldots,T/\dt~;
\label{Ahist}
\ee
note the similarity to \eqref{toepsol}.
As discussed in the introduction,
rather than an implicit solve of \eqref{Ahist} for each time step,
we prefer the following explicit scheme which
achieves the same order.

Consider the $k$th time step.
Firstly a ``predictor'' $\bmu^{(0)}$ is generated via a fixed order-$m$ extrapolation rule
in time applied to 
the density vector,
\be
\bmu^{(0)} \;=\; \sum_{r=1}^{m} c_r \bmu^{k-r}~,
\label{extrap}
\ee
where we choose $m=2q$ to match the D-spline order.
Here the $c_r$ are simply the values of the Lagrange polynomials associated with the time nodes
$t_{k-r}$ evaluated at $t_k$. 
Then $\tilde{\mbf{g}}^k$ is evaluated according to the right-hand side of
\eqref{Ahist}; this is the most expensive task.
Finally $n_c$ ``corrector'' steps are performed on $\bmu$, each of which is
a Jacobi iteration with shift $d$ ,as follows.
The system matrix is split into
the diagonal matrix $B$ and matrix $\tilde A^0$, defined by
\be
B_{jj} := \sfrac{1}{2} + A^0_{jj} + d~,
\qquad
\tilde A^0 := \sfrac{1}{2}I + A^0 - B~.
\label{B}
\ee
The linear system \eqref{Ahist} to be solved is then
$(\tilde A^0 + B)\bmu^k = \tilde{\mbf{g}}^k$.
Writing the $\alpha$th iterate for this solution as $\bmu^{(\alpha)} = \{\mu^{(\alpha)}_j\}_{j=1}^N$,
and initializing with \eqref{extrap},
the Jacobi iteration is
\be
\mu^{(\alpha+1)}_j \; = \; \frac{\tilde g_j^k - (\tilde A^0 \bmu^{(\alpha)})_j}{B_{jj}}
~,
\qquad j=1,\dots, N ~,
  \qquad \alpha = 0,\dots,n_c-1~.
\label{jacobi}
\ee
Once may interpret each iteration as decrementing $\mu_j$
by the $j$th component of its residual divided by $B_{jj}$.
This iteration, if it converges, converges to the exact (implicit) solution.
However, to make an explicit scheme we fix $n_c$, independent of $N$,
so that
the approximate solution to the time step is $\bmu^k = \bmu^{(n_c)}$.

This completes the description of the entire scheme for evolving
the density. The cost per time-step is $\bigO(N^2)$, and thus the total
$\bigO(N^2\cdot T/\dt) =\bigO(\dx^{-4}\dt^{-1})$. 
The numerical wave equation solution 
is then evaluated
with cost $\bigO(N)$
at any desired time and (not close) exterior target point
using \eqref{solrep} with the quadrature of Sec.~\ref{s:offsurf}.

\begin{rmk}
We chose a good diagonal shift empirically (by examining extremal eigenvalues)
as $d=-0.25$:
this vastly increases the corrector convergence rate in the
case of $\dt \gg \dx$, yet does no harm in other situations.
With this shift, $n_c=8$ was found to be sufficiently large to give
stability and errors similar to that of a full implicit solution.
Since $n_c \ll n$, the total corrector cost remains negligible
compared to that of evaluating $\tilde{\mbf{g}}^k$.
\label{r:shift}
\end{rmk}




\section{Numerical experiments}\label{NumExp} 

In this section we test the convergence of the full time-dependent
BIE scheme for exterior wave equation BVPs for the torus and cruller surfaces.
We set $n_r$ according to Remark~\ref{r:nr}, $n_\varphi = 2n_r$,
and $n_c$ and $d$ according to Remark~\ref{r:shift}.

\begin{figure}  
  \includegraphics[width=0.47\textwidth]{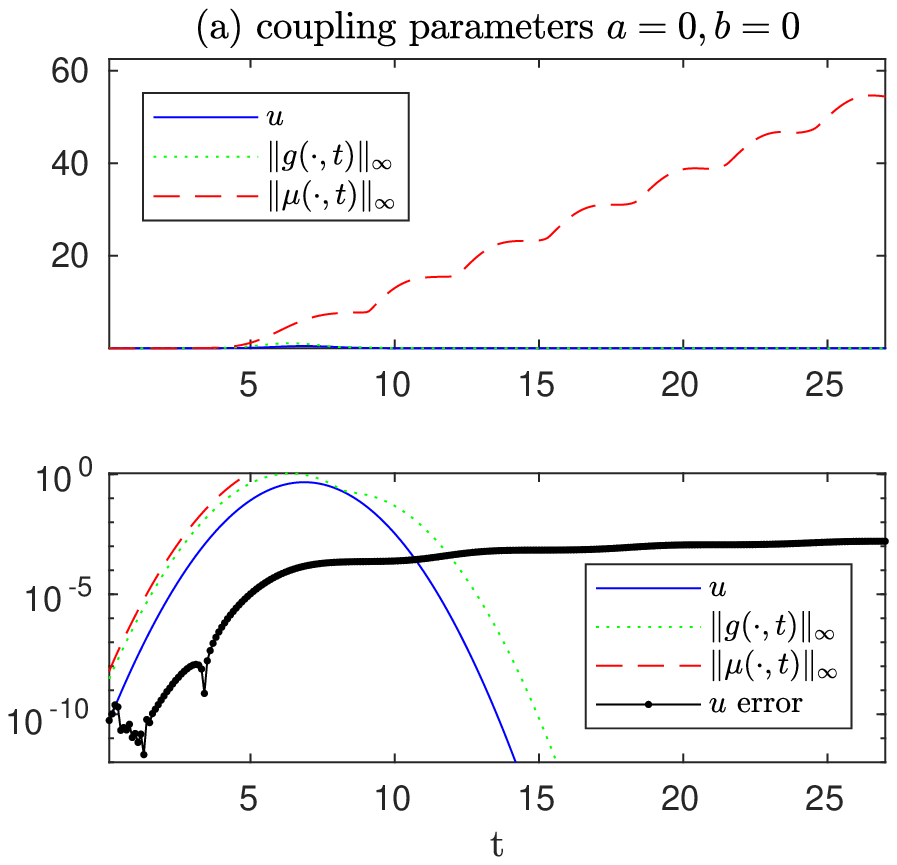}\qquad
  \includegraphics[width=0.47\textwidth]{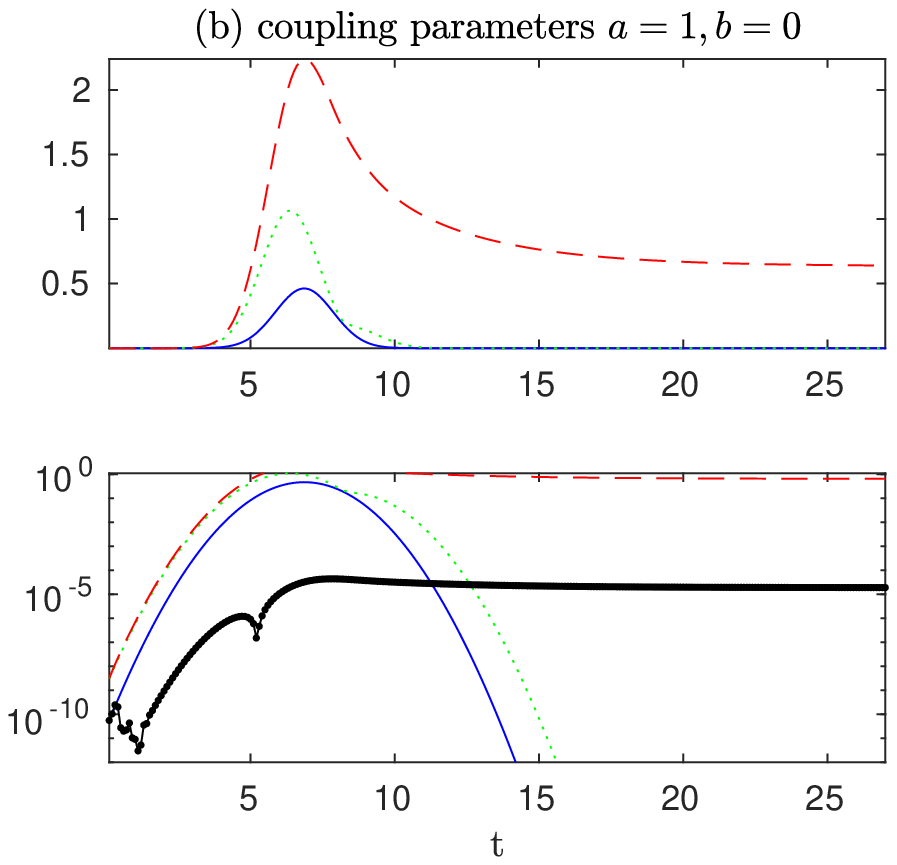}
  \\
  \\
  \includegraphics[width=0.47\textwidth]{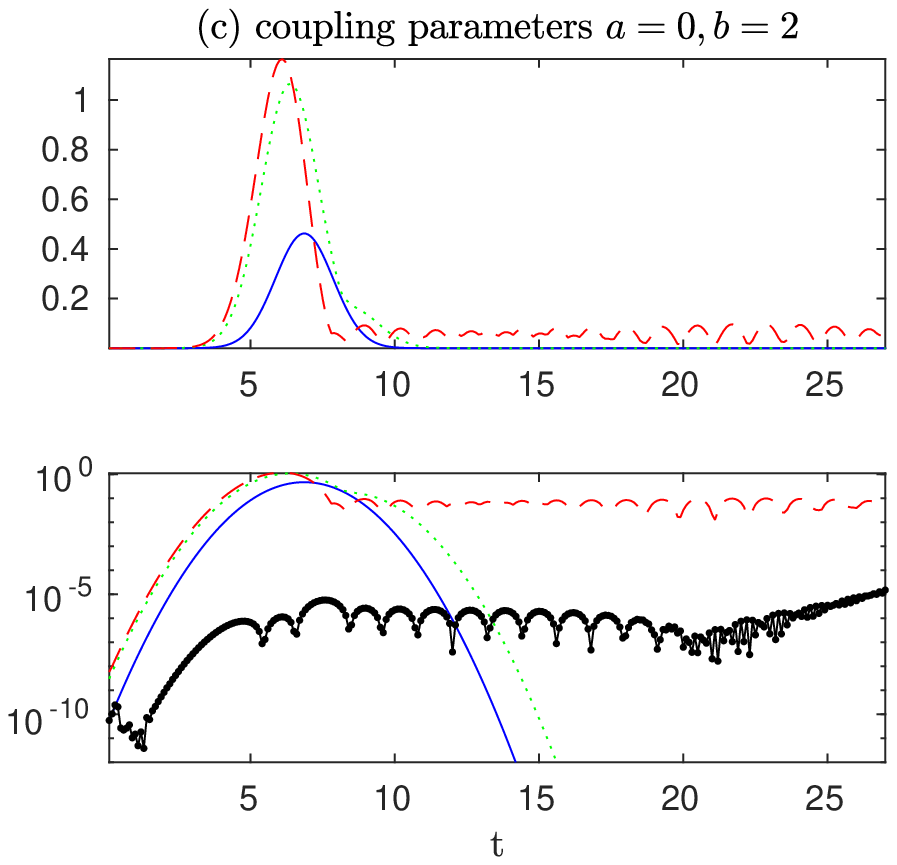}\qquad
  \includegraphics[width=0.47\textwidth]{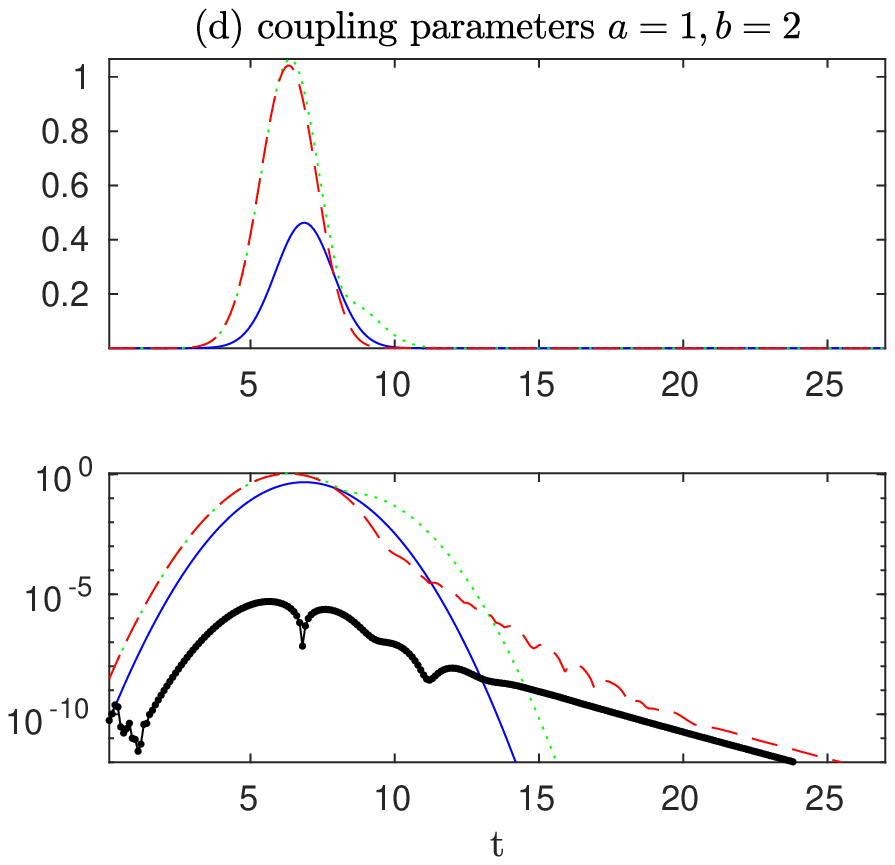}
  \caption{Time evolution of solution $u(x,t)$ (blue line)
    measured at the target
    point $x = (1.3,0.1,0.8)$ (shown in Fig.~\ref{f:spaceconv}(a)),
    the error in $u(x,t)$ (thick black line),
    the maximum over $\Gamma$ of the density $\mu$ (dashed red),
    and the maximum of the
    right hand side $g$ (dotted green), for a BVP with known solution in
    the exterior of the torus domain of Fig.~\ref{f:geom}(a).
    Four choices of coupling parameters $a$ and $b$ in the Volterra BIE scheme
    \eqref{cftdie} are shown, with $p=6$, $2q=4$,
    $N=1944$ ($\dx=0.108$) and $\dt=0.1$.
    Both linear (upper) and logarithmic (lower) vertical axes are shown.
    The true exterior solution is the retarded potential \eqref{usrc}
    from an interior monopole
    source $x_0 = (0.9,-0.2,0.1)$, emitting the temporal Gaussian pulse
    $T(t) = 5e^{-(t-6)^2/2}$.
  }
  \label{f:toggleab}
\end{figure}

\subsection{Effect of coupling parameters \texorpdfstring{$a$ and $b$}{{\it a} and {\it b}}}
\label{s:ab}

Choosing the torus surface, we pick an intermediate spatial order $p=6$,
with a low-resolution discretization of $\nph=9$ by $\nth=6$ panels,
thus $N=1944$, or $\dx=0.108$. For this experiment we
fix a timestep $\dt=0.1$,
with D-spline order $2q=p-2$ to match the spatial order.
We solve the exterior BVP with data $g$ deriving from the unit-magnitude
known solution given in the caption of Fig.~\ref{f:toggleab},
which, being a Gaussian pulse, dies away rapidly for $t>10$.
The four panel pairs of Fig.~\ref{f:toggleab} contrast the
resulting behavior of the norm of the density, and of the pointwise error
in $u(x,t)$,
for the four combinations of constant coupling parameters
$a \in \{0,1\}$ and $b\in\{0,2\}$.
The latter value $b=2$ is chosen as the maximum principal curvature of the torus.
These weights $a$ and $b$ are
motivated in the introduction, summarizing \cite{TDIEstab}.
Each combination results in a different behavior:
\bi
\item[(a)] $a=b=0$:
  Once the pulse has passed, the
  density $\mu$ grows asymptotically linearly in time.
  This secular growth 
  is associated with a zero-frequency Neumann resonance \cite{TDIEstab}.
  Accurate evaluation of $u$ for long times at any target point
  is thus impossible, due to growing catastrophic cancellation (indeed,
  a growing error in $u$ is visible).
\item[(b)] $a=1$, $b=0$:
  The situation is better than (a), with the size of $\mu$
  peaking then tending to a positive constant. The spatial function $\mu(x,t)$
  (not shown) tends to a constant on $\Gamma$; this is consistent
  with the decay of $u$ in \eqref{solrep} since ${\mathcal D}$ acting
  on a constant vanishes in $\Omega$.
\item[(c)] $a=0$, $b=2$: The situation is similar to (b), except that the
  density continues to oscillate with smaller constant amplitude around 0.1.
  At late times a weak instability (exponential growth in error) is seen.
\item[(d)] $a=1$, $b=2$:
  In this case only, density decays exponentially after the peak of the
  pulse. The decay rate transitions to a slower exponential rate
  (for $t>19$) once $\mu$ has dropped to the typical size of the error.
\ei
Thus the story for (a), (b), and (d) is exactly the same as found
previously for the sphere
\cite{TDIEstab} (case (c) was not tested in that work).
Only the last case leads to a density which dies
exponentially, thus the possibility of
potential evaluation with high relative accuracy.
The last case also leads to the smallest $\mu$ values and lowest errors in $u$:
the maximum $u$ error is $5 \times 10^{-6}$ (compare (b) for which it
is $4\times 10^{-5}$).


On a laptop with i7-7700HQ CPU, running a MATLAB implementation
(using Parallel Toolbox with 8 threads, calling
single-threaded Fortran90 for D-spline evaluation as in Sec.~\ref{TempD}),
the above calculation takes 16 seconds for the parallel
assembly of the sparse quadrature matrices (requiring around 10 GB),
then 7 s for each run (i.e.\ 0.025 s per time step),
which is dominated by the single-threaded sparse matrix-vector products
which evaluate $\tilde{\mbf{g}}^k$ in the right-hand side of \eqref{Ahist}.

We remark that the exterior solution in this case decays much more rapidly than a typical solution of a scattering problem, which would be no faster than exponential, as determined by the scattering resonances. Therefore here one cannot demand a density decay rate matching that of the solution. In the more realistic scattering 
problem considered below, the decay rates are better matched. 

\begin{figure}  
  \mbox{\hspace{-8ex}\includegraphics[width=1.2\textwidth]{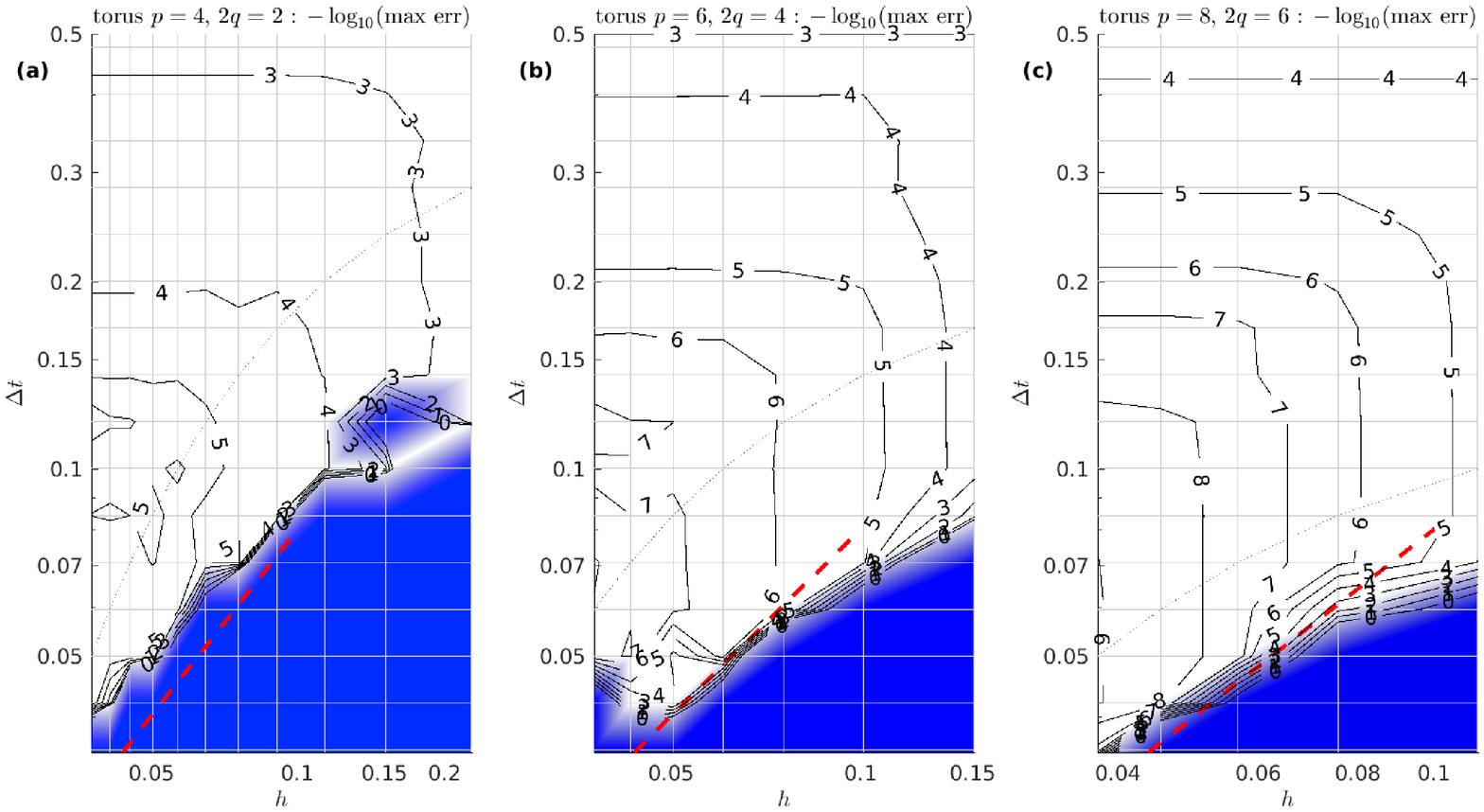}}
  \\
  \vspace{-5ex}
  \\
  \mbox{\hspace{-8ex}\includegraphics[width=1.2\textwidth]{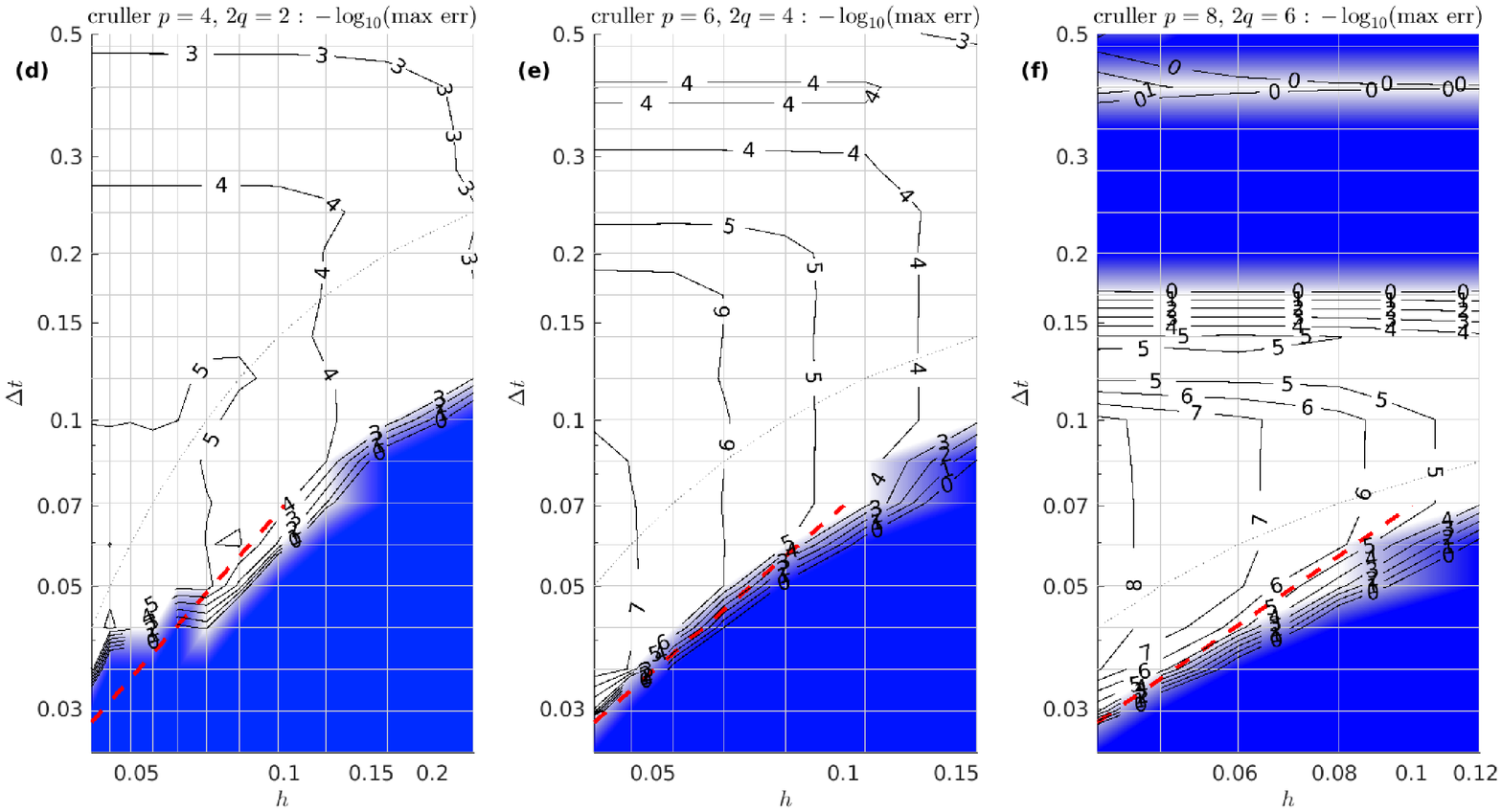}}
  \vspace{-6ex}
  \caption{
    Convergence with respect to resolution $\dx$ (see \eqref{dxdef})
    and time-step $\dt$
    of the maximum error in $u(x,t)$ at the target $x=(1.3,0.1,0.8)$, for the
    same exterior
    BVP as in Fig.~\ref{f:toggleab}. Couplings are $a=1$, $b=2$.
    Scheme orders increase along each row.
    Top: torus; bottom: cruller (see Fig.~\ref{f:geom}).
    Each gray grid shows $(\dx,\dt)$ pairs tested.
    Contours show the number of correct digits,
    i.e.\ $-\log_{10} \max_{t} |\tilde u(x,t)- u(x,t)|$, where $\tilde u$ is
    the numerical solution with parameters $(\dx,\dt)$.
    Blue shaded regions show where the error exceeds 1, indicating
    failure and, usually, instability.
    A rough stability boundary $\dt = c_\tbox{iCFL}\dx$
    (dashed red line) is shown, with slope
    $c_\tbox{iCFL}=0.8$ for the torus, $0.7$ for the cruller.
    ``$\dx$-dominated curves'' used to measure $\dx$-convergence
    are shown (dotted gray).}
  \label{f:stab}
\end{figure}

\begin{figure}  
  \includegraphics[width=1.0\textwidth]{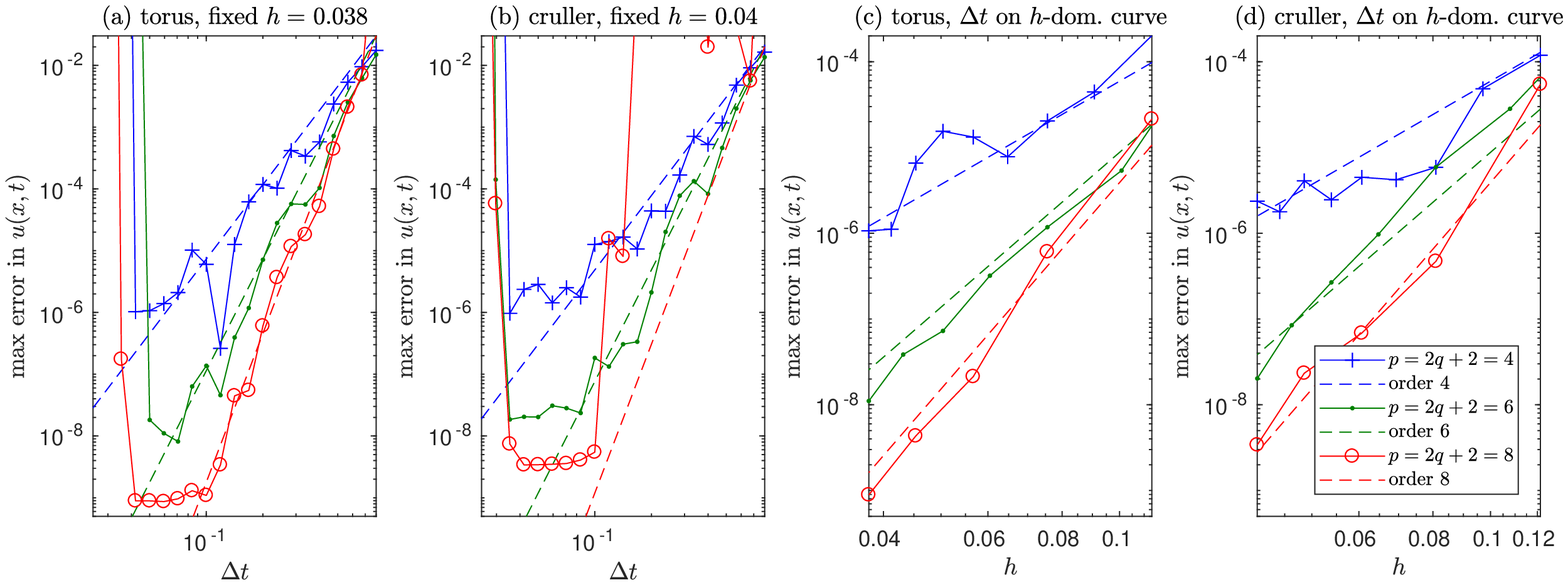}
  \vspace{-5ex}
  \caption{Convergence of maximum error in $u(x,t)$
    with respect to time-step $\dt$
    (panels a--b) and resolution $\dx$ (panels c--d),
    for the same exterior BVP as in Fig.~\ref{f:toggleab}.
    See legend of panel (d) for the three different orders tested
    (solid lines) and the expected orders (dashed lines).
    (a) and (c) are for the torus, and (b) and (d) for the cruller.
    For (c) and (d), in order to avoid the unstable small-$\dt$ region,
    $\dt$ is scaled with $h$ along the grey dotted curves of
    Fig.~\ref{f:stab}.}
  \label{f:orders}
\end{figure}

\subsection{Stability and convergence}
\label{s:stab}

Here we report on tests of accuracy and stability covering the $(\dx,\dt)$
plane, for both torus and cruller domains.
The same Gaussian pulse BVP as in Sec.~\ref{s:ab} was used, and the same
target point.
The coupling parameters were fixed at $a=1$ and $b=2$ as motivated above.
The temporal order $2q+2$ was set to match the spatial panel order $p$,
with orders 4, 6, and 8 tested.

Fig.~\ref{f:stab} shows the resulting numbers of accurate digits.
Note that the peak in $u(x,\cdot)$ at the target $x$ is of order 1;
thus panel (c) shows over 9-digit relative accuracy for $p=8$ at the minimum $\dx$.
The main results are:
\bi
\item In every case there is an unstable region (lower right, shaded in blue)
  in the approximate form of an inverse-CFL condition \eqref{iCFL}
  for some $\bigO(1)$ constant $c_\tbox{iCFL}$
  that appears to be at most weakly dependent upon the order
  and upon the domain.
  Dashed red lines show a boundary of the form \eqref{iCFL}, to guide the
  eye; the contours and shading show some minor variations from this
  form.
\item For the cruller at high orders, there is an additional stability
  constraint that $\dt$ not exceed a certain ($h$-independent) value $\dt_\tbox{max}$. For $p=8$ it appears that $\dt_\tbox{max} \approx 0.15$.
  For the torus, the scheme is stable for $\dt$ as large as 0.8 (a quarter of the domain diameter, and much bigger than the typical panel size).
\item For any fixed $\dt<\dt_\tbox{max}$,
  the scheme appears stable as $\dx\to 0$.
\item When stable,
  convergence with the expected high orders is seen, with respect
  to both $\dt$ and $\dx$. The contours of constant error bend sharply
  from horizontal to vertical in each $(\dx,\dt)$ plane, indicating the
  transition from $\dt$-dominated error (upper left region) to
  $\dx$-dominated (lower right). Fig.~\ref{f:orders} quantifies the
  convergence rates: panels (a) and (b) show $\dt$-convergence at the
  fixed smallest $\dx$, matching the expected temporal order $2q+2$
  in all stable regions.
  Because of the inverse-CFL condition, one cannot fix a small $\dt$ to
  test the $\dx$-convergence; instead we varied $\dt$ with $\dx$ along
  ``$\dx$-dominated'' curves, shown (grey dotted) in Fig.~\ref{f:stab},
  which remains stable and in a region where contours are vertical.
  Panels (c) and (d) show that the resulting $\dx$-convergence is, barring
  some variation, consistent with $\bigO(\dx^p)$, as expected for
  the interpolation order (Sec.~\ref{s:predcorr}).
  The $\bigO(\dx^{2p})$ spatial quadrature order is no longer seen.
\ei

\begin{rmk}[Other assessments of stability]
  Stability was also assessed by extracting the growth/decay rate
  of $\|\mu(\cdot,t)\|_\infty$ near the end of the time interval, $T=50$.
  The results were consistent with the blue shading of Fig.~\ref{f:stab}
  extracted simply from the error, so are not shown.
  More sophisticated attempts to assess stability did not prove
  useful.
  For instance, setting $\mbf{g}^k=0$ in \eqref{Ahist}, the entire
  predictor-corrector scheme, followed by a backwards shift by 1 in $k$,
  can be viewed as a huge sparse square matrix acting on the density
  history vector $\{\bmu^k, \bmu^{k-1},\dots,\bmu^{k-n}\} \in \RR^{nN}$.
  We estimated extremal eigenvalues of this matrix via
  ARPACK \cite{arpack} ({\tt eigs} in MATLAB),
  but convergence was extremely slow, often
  taking much more time than the $T/\dt$ time-steps needed
  (we speculate that this is due to its block-companion matrix structure).
\end{rmk}

Collecting the data in Fig.~\ref{f:stab} required 966 solution runs,
with between $63$ and $2000$ time-steps per run,
and $N$ varying from $384$ to $13824$.
This required 11 CPU-days on a server with two 14-core
Intel Xeon 2.4GHz 2680v4 CPUs and 512 GB of RAM.
Most of the time is taken
by time-stepping (single-threaded sparse matrix-vector products).
Because of our direct $\bigO(\dx^{-4})$ implementation,
run time is dominated by the one or two smallest $\dx$ values.
The RAM cost of $\bigO(\dx^{-4})$, dominated by forming and storing the set $\{A^r\}_{r=0}^n$, limited the minimum $\dx$ that could be tested;
at our minimum $\dx$, 260 GB of RAM was used.

\begin{figure}  
  \mbox{\includegraphics[width=0.5\textwidth]{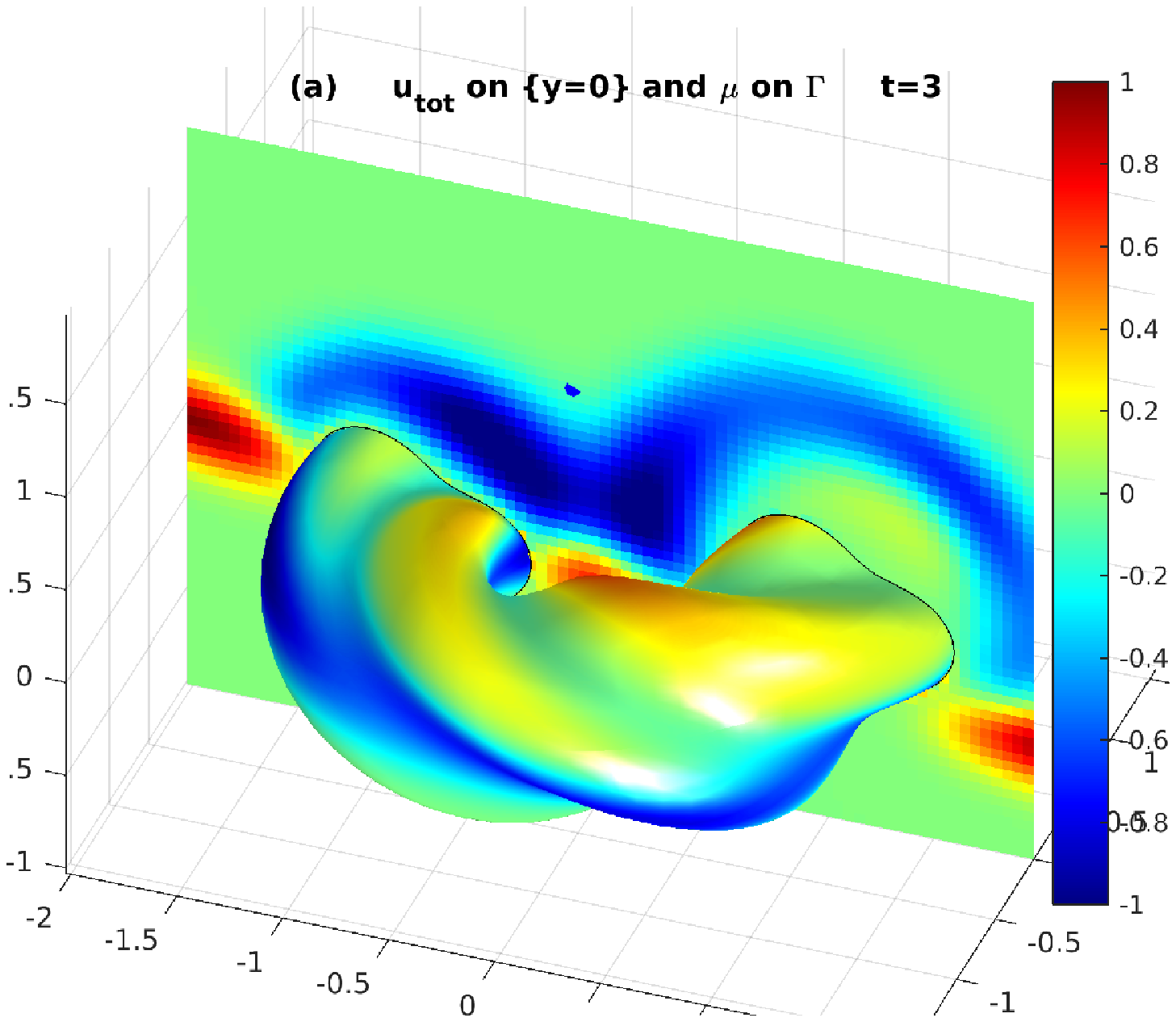}
  \includegraphics[width=0.5\textwidth]{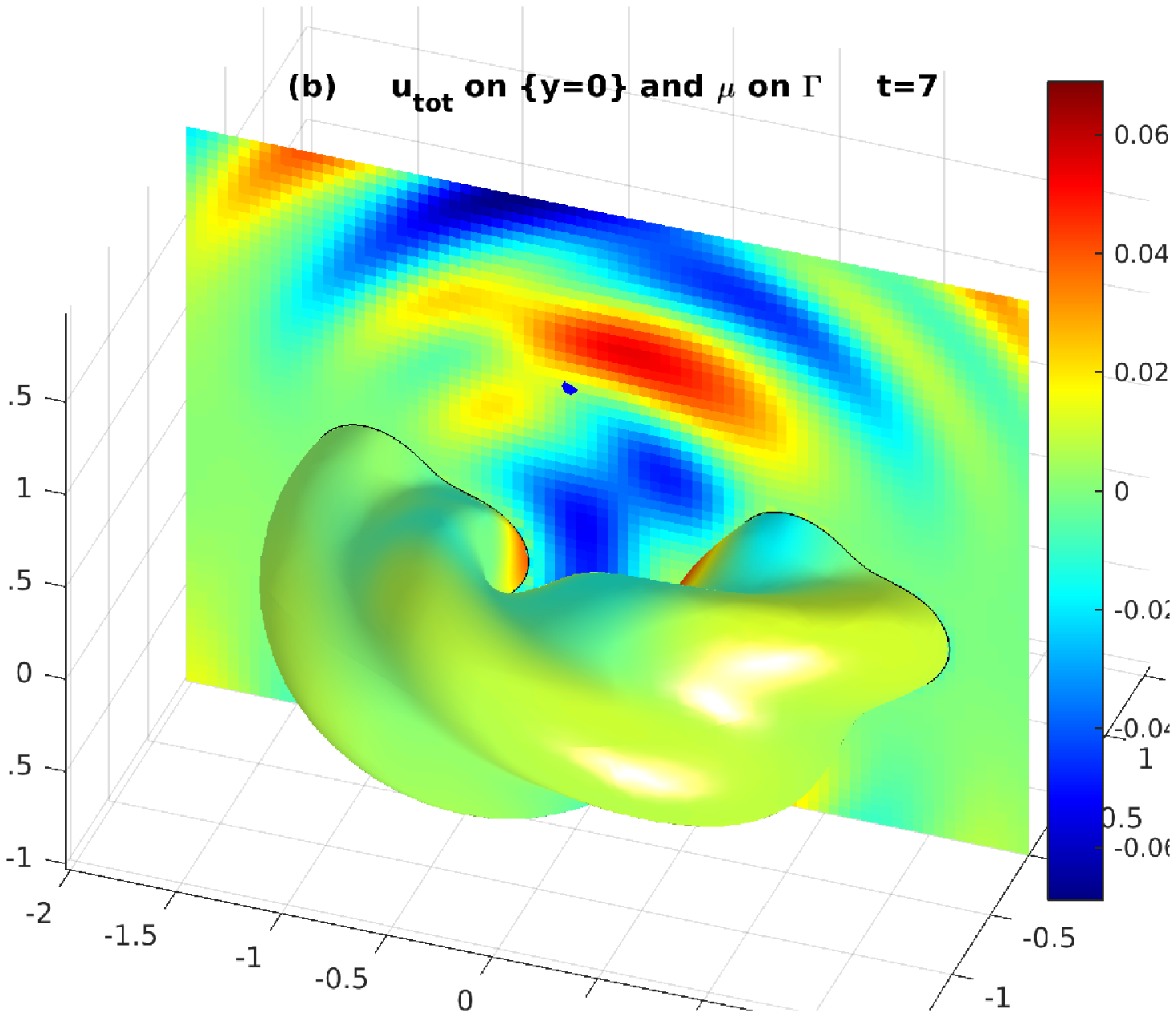}}
  \vspace{-1ex}
  \\
  {\centering
  \includegraphics[width=0.48\textwidth]{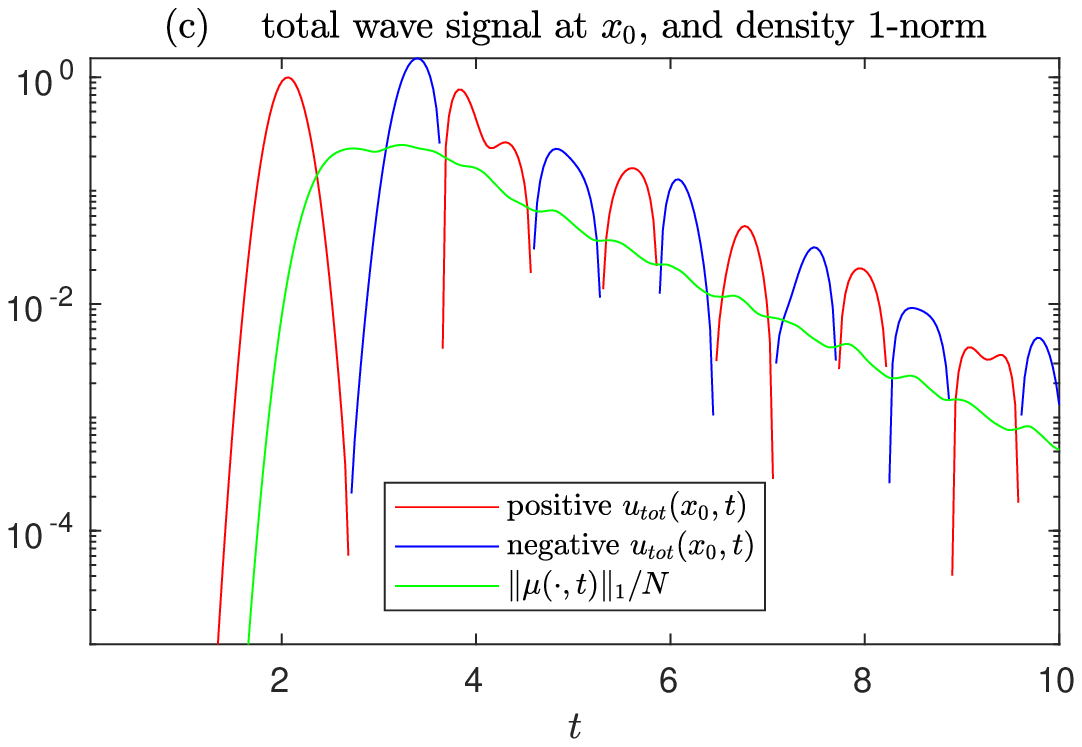}
  \;
  \includegraphics[width=0.5\textwidth]{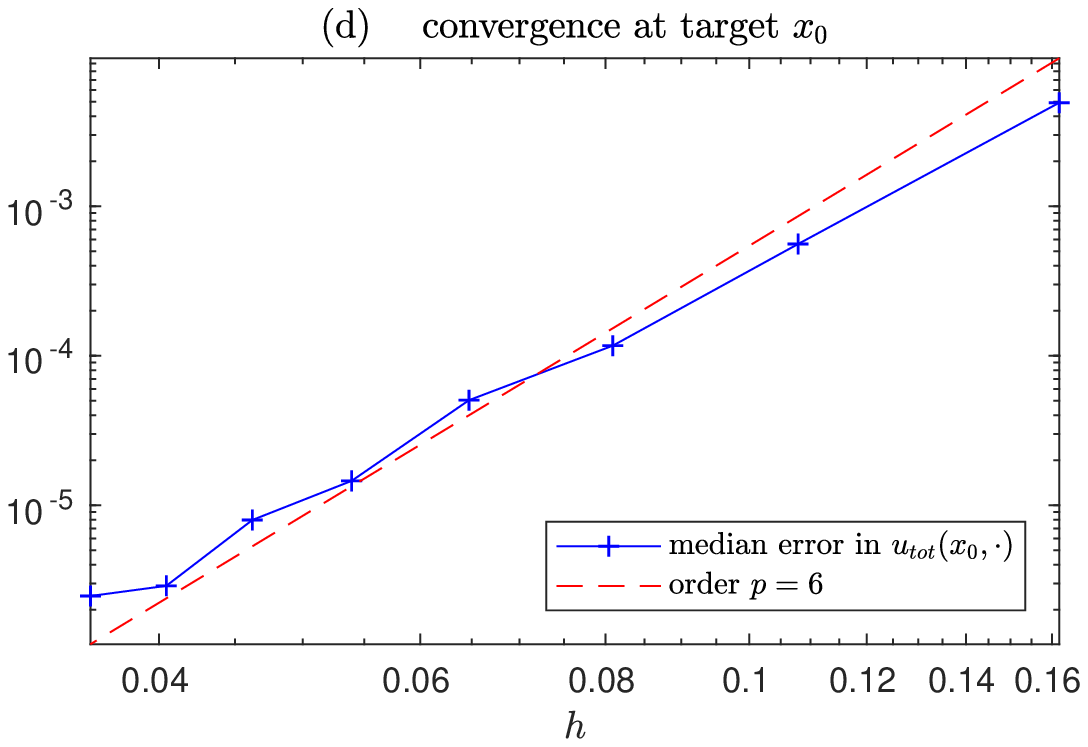}
  }
  \vspace{-2ex}
  \caption{Scattering of an incident plane wave from a sound-soft cruller.
    The unnormalized
    incident wave direction is $\mbf{d} = (-0.2,0.1,-1)$, i.e.\ the wave
    comes down from above, with temporal
    pulse $T(t) = e^{-(t-3)^2/2\sigma^2}$ where $\sigma=0.15$.
    (a) and (b) show the full wave $u_\tbox{tot} = u_\tbox{inc}+u$ evaluated on an array
    of 4941 points on the slice $\{y=0\}$, and $\mu(\cdot,t)$ on $\Gamma$.
    In (a) the incident wave has just hit the obstacle;
    in (b) only remnant decaying radiation is visible (note the color scale).
    (c) shows the time dependence of $u_\tbox{tot}$ at the point $x_0$ shown
    as a black dot in (a) and (b); the peak at $t\approx 2$ is the incident
    wave before collision, and the reflected signal is at $t>3$.
    (d) shows median difference in the solution $u_\tbox{tot}(x_0,t)$
    from its converged values,
    when interpolated to a $t$ grid in $[1,9]$, 
    as a function of resolution $\dx$. We fix $\dt = 1.0 \dx$.
    }
  \label{f:scatt}
\end{figure}

\subsection{Plane wave pulse acoustic scattering example}
\label{s:scat}

Finally we demonstrate the solver for the Dirichlet (sound-soft) scattering
application, using the cruller domain of Fig.~\ref{f:geom}(b).
The boundary data is $g = -u_\tbox{inc}$
on $\Gamma$, for the plane wave
$$
u_\tbox{inc}(x,t) \; = \; T(t - \hat{\mbf{d}}\cdot x)~,
\qquad \hat{\mbf{d}} := \mbf{d} / |\mbf{d}|~,
$$
where $\mbf{d}$ defines its direction and
$T(t)$ its signal function.
$u_\tbox{inc}$ is a solution to the wave equation in
$\RR^3 \times \RR$. Fig.~\ref{f:scatt}(a--b) shows
two snapshots of
the resulting physical (full) wave $u_\tbox{tot} = u_\tbox{inc}+u$,
which vanishes on $\Gamma$ for all time,
with $\mbf{d}$ and $T(t)$ as stated in the
figure caption.
Note that the pulse $T(t)$ is about 7 times narrower than the pulse
used in previous tests.
Its full-width half-maximum is about 0.35,
i.e.\ around 1/9 of the domain diameter.

Our experiment is done at spatial
order $p=6$, with $2q=4$ to match in temporal order,
and $n_r = 2p+4$, $n_\varphi =2n_r$ as in Remark~\ref{r:nr}.
Fig.~\ref{f:scatt}(d) shows the pointwise convergence with respect to
$\dx$, for $\dt = 1.0\dx$.
The median error is shown estimated by comparison against the smallest
$\dx$ case ($\nph=30$).
It is quite consistent with 6th order convergence.
The smallest-$\dx$ point shown,
with median estimated error $2.5\times 10^{-6}$,
used $\nph = 27$ by $\nth=18$ panels, thus $N = 17496$ spatial nodes.
This took 44 minutes on the Xeon server mentioned above,
using largely one thread.
Of this, 23 minutes was sparse matrix filling (needing
around two billion nonzero elements in $\{A^r\}$ matrices),
and 21 minutes for the 279 time steps. 290 GB of RAM was used.
%

We note that if only 3-digit accuracy is needed, 
$n_\phi = 9$ ($\dx \approx 0.11$) is sufficient,
requiring only 1 minute, and 11 GB of RAM, on the
laptop mentioned in Sec.~\ref{s:ab}.

An overall exponential decay of the solution is indicated by
Fig.~\ref{f:scatt}(c); notice that the signal also oscillates as it decays.
Although domains without trapped rays have long been known to have exponential
decay of the wave equation solution \cite{morawetz77},
the cruller has trapped rays---for example it has one in the plane $z=0$ reflecting off the five symmetric bumps encircling the hole.
Assuming the hard-walled case is similar to the case of smooth potentials,
recent analysis predicts exponential decay, in the high-frequency limit, at a rate determined by the smallest distance of scattering resonances from the real axis, which is controlled by the weakest trapped ray Lyapunov exponent \cite{NZ09} \cite[Ch.~6]{DZ19}.
The figure also plots the decay of the density 1-norm: it appears to
be exponential at the same rate as the solution,
which is thus the optimal rate. This indicates an absence of catastrophic
cancellation in the solution representation, a major advantage of our formulation.

\section{Conclusions and future work}\label{Conclude} 

We have implemented a recently-proposed
combined-field time-domain integral equation
formulation \cite{TDIEstab}
to solve to high-order accuracy
the Dirichlet BVP for the scalar wave equation
in a general smooth exterior domain in $\RR^3$.
This showcases a method previously only studied for the sphere (exploiting
separation of variables \cite{TDIEstab}).
We provide evidence that in relatively general domains, as for the
sphere, both of the coupling weights $a$ and $b$ must be positive to prevent
long-lived resonances that cause catastrophic cancellation
with conventional formulations.
The retarded potentials are discretized in space
using high-order quadratures for weakly-singular kernels,
using new difference-spline temporal interpolants.
The timestepping is explicit, via a predictor-corrector scheme.
We verified the expected high order accuracy,
explored stability in the $(\dx,\dt)$ plane---which indicates an inverse CFL condition as found in the modal sphere case---%
and benchmarked a plane wave scattering example in which the exponential density
decay rate appears to match that of the solution.
In contrast to direct discretization methods,
time-domain integral formulations need not have any CFL condition
(upper bound on $\dt/\dx$), and we do not observe one.

The implementation presented is direct, 
since it needs $\bigO(N^2)$ memory and effort per timestep,
and assumes a structured grid of quad patches that is appropriate only
for torus-like surfaces.
In this setting, because of our high orders up to 8th,
we achieve typically 5--9 digits of accuracy.
Thus, for low-frequency-content data, needing $N \lesssim 10^4$,
the scheme is quite useful for high-accuracy solutions on smooth surfaces.

The work suggests several directions for improvement.
\bi
\item
  Handling more general shapes would require
  unstructured, adaptive triangles and/or quad patches.
  One route to handle the kernel singularity is then to borrow from 
  on-surface quadrature schemes for harmonic layer potentials 
  of this same singularity class
  \cite{brunoFMM,ying06,bremer3d}.
\item
  To get high accuracy evaluation near to the surface, a special quadrature
  scheme would be needed, building on a variety of existing schemes for the 
  Laplace equation (see, for example, \cite{qbx}).
\item
  In order to effectively address large-scale problems, our marching scheme
  should be coupled with the fast algorithms mentioned in the introduction 
\cite{Bleszynski,PWFTD,Chewfastbook,MengBoagMichielssen,Shanker2003,Yilmaz2004}.
\item
  Extensions of our formulation and discretization to electromagnetic scattering problems using Debye sources, as discussed in \ci{DebyeSource,Sphere_EM}, 
seem possible.
\item
  One may be able to optimize stability, or density decay rate, by a better
  choice of positive coupling parameters, or functions, $a$ and $b$.
\item
  In our studies we tied the temporal order $2q+2$ to match the
  spatial order $p$, and observed a somewhat restrictive
  $\dt_\tbox{max}$ for the cruller domain with $p=8$.
  Independently varying $2q$ and $p$ would expose whether this
  constraint is tied to the temporal or spatial high order.
\item
  There are many opportunities for rigorous analysis of the discretization
  scheme.
\ei

Finally, we note that the inverse-CFL condition
may pose a problem for spatially-adaptive quadratures
in the low-frequency regime,
since the large panels would place a lower bound on $\dt$ that would
be inaccurate for the small features.
It is possible that a modified temporal interpolant could remove this
condition, or that locally-adaptive timestepping could circumvent it.
In particular, since the temporal interpolants are constructed independently for each point on the spatial grid, it is possible to use different time grids wherever needed. However the stability of such a scheme requires further study.

\section*{Acknowledgments}
We are grateful for useful discussions with Charlie Epstein, Maciej Zworski, and Andrew Hassell.

\bibliographystyle{plain}
\bibliography{refs}

\end{document}